\documentclass[amstex,12pt,thmsa,sw20lart]{article}
\usepackage{amsfonts, amsmath,amsthm,amssymb,color}

\input amssym.def
\input amssym.tex

\begin{document}

\date{}
\title{{\bf Pointwise products of  some Banach \\
function spaces and factorization}}

\author{Pawe{\l } Kolwicz\thanks{Research partially supported by the State Committee 
for Scientific Research, Poland, Grant N N201 362236}\,, \ Karol Le\'{s}nik$^\ast$ and
Lech Maligranda}
\date{}
\maketitle

\vspace{-7mm}

\begin{abstract}

\noindent {\footnotesize The well-known factorization theorem of Lozanovski{\u \i} may 
be written in the form $L^{1}\equiv E\odot E^{\prime }$, where $\odot $ means the
pointwise product of Banach ideal spaces. A natural generalization of this problem 
would be the question when one can factorize $F$ through $E$, i.e., when 
$F\equiv E\odot M(E, F) \,$, where $M(E, F) $ is the space of pointwise multipliers from 
$E$ to $F$. Properties of $M(E, F) $ were investigated in our earlier paper [KLM12] and 
here we collect and prove some properties of the construction $E\odot F$. The formulas 
for pointwise product of Calder\'{o}n-Lozanovski{\u \i}  $E_{\varphi}$ spaces, Lorentz spaces 
and Marcinkiewicz spaces are proved. These results are then used to prove factorization 
theorems for these spaces. Finally,  it is proved in Theorem 11 that under some natural 
assumptions, a rearrangement invariant Banach function space may be factorized through 
Marcinkiewicz space.  }
\end{abstract}

\renewcommand{\thefootnote}{\fnsymbol{footnote}}


\vspace{-3mm}
\begin{center}
\textbf{1. Introduction and preliminaries}
\end{center}

The well-known factorization theorem of Lozanovski{\u \i} says that for any $\varepsilon > 0$ 
each $z\in L^{1}$ can be factorized by $x \in E$ and $y \in E^{\prime }$ in such a
way that 
\begin{equation*}
z = x y ~~{\rm and} ~~ \| x \|_{E}\, \| y \|_{E^{\prime }}\leq (1+\varepsilon ) \| z \|_{L^{1}}.
\end{equation*} 
Moreover, if $E$ has the Fatou property we may take $\varepsilon =0$ in the above inequality. This
theorem can be written in the form $L^{1}\equiv E\odot E^{\prime }$, where 
\begin{equation*}
E\odot F = \left\{ x\cdot y: x\in E\text{ and }y\in F\right\} .
\end{equation*}
Then natural question arises: when is it possible to factorize $F$ through $E$, i.e., when 
\begin{equation}
F\equiv E\odot M\left( E, F\right) ?  \label{factoriz}
\end{equation}
Here $M\left( E, F\right) $ denotes the space of multipliers defined as 
\begin{equation*}
M\left( E, F\right) =\left\{ x\in L^{0}: x y \in F\text{ for each }y\in E\right\}
\end{equation*}
with the operator norm 
\begin{equation*}
\left\| x\right\|_{M\left( E, F\right) }=\sup_{\left \| y\right\|_{E}=1}\left\| xy\right\|_{F}.
\end{equation*}

The space of multipliers between function spaces was investigated by many
authors, see for example [Ru79], [MP89] and [AZ90] (see also  [Za66],  [ZR67], 
[Cr72], [Ma74], [AS76], [Ma89], [Ra92], [Na95], [DR00], [CN03], [CDS08], [MN10], 
[Sc10] and  [KLM12]). In this paper we are going to investigate general properties 
of the product construction $E\odot F$ and calculate the product space $E\odot F$ for 
Calder\'{o}n-Lozanovski{\u \i}, Lorentz and Marcinkiewicz spaces. This product space 
was of interest in [An60], [Wa63], [ON65], [ZR67], [Da74], [Ru79], [Ma89], [RR91], [Ra92], 
[BL93], [DMM03], [AM09], [KM10] and [Sc10]. The results on product construction will be 
used to give answers to the factorization question (\ref{factoriz}) in these special spaces.

Let $(\Omega ,\Sigma ,\mu )$ be a complete $\sigma $-finite measure space
and $L^{0}=L^{0}(\Omega )$ be the space of all classes of $\mu $-measurable
real-valued functions defined on $\Omega $. A (quasi-) Banach space 
$E=\left( E,\| \cdot \| _{E}\right) $ is said to be a ({\it quasi-}) {\it Banach ideal space} 
on $\Omega $ if $E$ is a linear subspace of $L^{0}(\Omega )$ and satisfies the 
so-called {\it ideal property}, which means that
if $y\in E,x\in L^{0}$ and $|x(t)|\leq |y(t)|$ for $\mu $-almost all $t\in
\Omega $, then $x\in E$ and $\| x\|_{E}\leq \| y\|_{E}$. We
will also assume that a ({\it quasi-}) {\it Banach ideal space} on $\Omega $ is
saturated, i.e. every $A\in \Sigma $ with $\mu (A)>0$ has a subset $B\in
\Sigma $ of finite positive measure for which $\chi _{B}\in E$. The last
statement is equivalent with the existence of a {\it weak unit}, i.e., an
element $x\in E$ such that $x( t) >0$ for each $t\in \Omega $
(see [KA77] and [Ma89]). If the measure space $(\Omega ,\Sigma,\mu )$ 
is non-atomic we should say about ({\it quasi-}) {\it Banach function
space}, if we replace the measure space $\left( \Omega ,\Sigma ,\mu \right) $
by the counting measure space $\left( \mathbb{N},2^{\mathbb{N}},m\right) ,$
then we say that $E$ is a ({\it quasi-}) {\it Banach sequence space} (denoted
by $e$). 

A point $x\in E$ is said to have {\it order continuous norm} (or to be
{\it order continuous element}) if for each sequence $(x_{n})\subset E$
satisfying $0\leq x_{n}\leq |x|$ and $x_{n}\rightarrow 0~\mu $-a.e. on 
$\Omega ,$ one has $\| x_{n} \|_{E}\rightarrow 0$. By $E_{a}$ we denote
the subspace of all order continuous elements of $E$. It is worth to notice
that in case of Banach ideal spaces on $\Omega $, $x\in E_{a}$ if and only
if $\| x\chi _{A_{n}}\|_{E}\downarrow 0$ for any sequence $\{A_{n}\}$
satisfying $A_{n} \searrow \emptyset$ (that is $A_{n}\supset A_{n+1}$ and 
$\mu (\bigcap_{n=1}^{\infty }A_{n})=0$). A Banach ideal space $E$ is called 
{\it order continuous} if every element of $E$ is order continuous, i.e., 
$E = E_{a}$.

A space $E$ has the {\it Fatou property} if $0\leq x_{n}\uparrow x\in
L^{0}$ with $x_{n}\in E$ and $\sup_{n\in \mathbb{N}}\Vert x_{n}\Vert
_{E}<\infty $ imply that $x\in E$ and $\| x_{n} \|_{E}\uparrow \|
x \|_{E}$.

We shall consider pointwise product of Calder\'{o}n-Lozanovski{\u \i} 
spaces $E_{\varphi }$ (which are generalizations of Orlicz spaces 
$L^{\varphi }$), but very useful will be also general Calder\'{o}n-Lozanovski{\u \i} 
construction $\rho ( E, F) $.

A function $\varphi: [0,\infty )\rightarrow \lbrack 0,\infty ]$ is called a 
{\it Young function} (or an {\it Orlicz function} if it is finite-valued)
if $\varphi $ is convex, non-decreasing with $\varphi (0)=0$; we assume also
that $\varphi $ is neither identically zero nor identically infinity on 
$(0,\infty )$ and $\lim_{u\rightarrow b_{\varphi }^{-}}\varphi (u)=\varphi
(b_{\varphi })$ if $b_{\varphi }<\infty $, where $b_{\varphi }=\sup \{u>0:\varphi (u)<\infty \}.$

Note that from the convexity of $\varphi $ and the equality $\varphi (0)=0$
it follows that $\lim_{u\rightarrow 0+}\varphi (u) \newline=\varphi (0)=0$.
Furthermore, from the convexity and $\varphi \not\equiv 0$ we obtain that 
$\lim_{u\rightarrow \infty }\varphi (u)=\infty $.

If we denote $a_{\varphi }=\sup \{u\geq 0:\varphi (u)=0\}$, 
then $0\leq a_{\varphi }\leq b_{\varphi }\leq \infty $ and $a_{\varphi
}<\infty ,~b_{\varphi }>0$, since a Young function is neither identically
zero nor identically infinity on $(0,\infty )$. The function $\varphi $ is
continuous and nondecreasing on $[0,b_{\varphi })$ and is strictly
increasing on $[a_{\varphi },b_{\varphi })$. If $a_{\varphi }=0$ then we
write $\varphi >0$, if $b_{\varphi }=\infty ,$ then $\varphi <\infty $. 
Each Young function $\varphi $ defines the function 
$\rho_{\varphi }:[0,\infty )\times \lbrack 0,\infty )\rightarrow \lbrack 0,\infty)$ 
in the following way
\begin{equation*}
\rho _{\varphi }(u,v)=\left\{ 
\begin{array}{cc}
v\varphi ^{-1}(\frac{u}{v}) & \text{if }u>0, \\ 
0 & \text{if }u=0,
\end{array}
\right.  \label{ro przez fi}
\end{equation*}
where $\varphi^{-1}$ denotes the right-continuous inverse of $\varphi$ and is defined by 
$$
\varphi ^{-1}(v) = \inf \{u\geq 0: \varphi (u)>v\}~\mathrm{for}~v\in \lbrack
0,\infty )~~\mathrm{with}~~\varphi ^{-1}(\infty )=\lim_{v\rightarrow \infty
}\varphi ^{-1}(v).
$$ 
If $\rho =\rho_{\varphi }$ and $E, F$ are Banach ideal spaces over the
same measure space $(\Omega, \Sigma, \mu)$, then the {\it Calder\'{o}n-Lozanovski{\u \i} 
space} $\rho (E, F)$ is defined as all $z \in L^0(\Omega)$ such that for some $x \in E, y \in F$ 
with $\| x\|_E \leq 1, \| y\|_F \leq 1$ and for some $\lambda > 0$ we have
\begin{equation*}
|z | \leq \lambda\, \rho(|x|, |y|) ~~ \mu-{\rm a.e. ~on} ~~ \Omega. 
\end{equation*}
The norm $\|z\|_{\rho} = \| z\|_{\rho(E, F)}$ of an element $z \in \rho(E, F)$ is defined as the 
infimum values of $\lambda$ for which the above inequality holds. It can be shown that
\begin{equation*}
\rho (E, F) = \left\{ z\in L^{0}(\Omega): |z| \leq \rho (x, y) ~ \mu-{\rm a.e. ~for ~some} ~ x \in E_+, 
~ y \in F_+ \right \}
\end{equation*}
with the norm 
\begin{equation}
\|z \|_{\rho \left( E, F\right) }=\inf \left\{ \max \left\{\| x \|_{E}, \| y \|_{F}\right\}: \text{\ }
| z| \leq \rho (x, y)\text{,}\ x \in E_+, y \in F_+ \right\} \text{.}  \label{def norm CL}
\end{equation}
The Calder\'{o}n-Lozanovski{\u \i} spaces, introduced by Calder\'{o}n in [Ca64] and 
developed by Lozanovski{\u \i} in [Lo71], [Lo73], [Lo78a] and [Lo78b], play crucial 
role in the theory of interpolation since such construction is interpolation functor for 
positive operators and under some additional assumptions on spaces $E, F$, like 
Fatou property or separability, also for all linear operators (see [Ov76], [Ov84], [KPS82], [Ma89]). 
If $\rho (u,v)=u^{\theta }v^{1-\theta} $ with $0<\theta <1$ we write $E^{\theta }F^{1-\theta }$ 
instead of $\rho (E,F)$ and these are Calder\'on spaces (cf. [Ca64], p. 122). Another important 
situation, investigated by Calder\'on (cf. [Ca64], p. 121) and independently by Lozanovski{\u \i} 
(cf. [Lo64, Theorem 2], [Lo65, Theorem 2]), appears when we put $F\equiv L^{\infty }$. 
In this case, in the definition of  the norm, it is enough to take $y=\chi_{\Omega }$ and then 
\begin{equation}
\| z \|_{\rho_{\varphi }\left( E, L^{\infty }\right) } = \inf
\left\{ \lambda > 0: \left \| \varphi  \left( | x |
/\lambda \right) \right \|_{E}\leq 1\right\}.  \label{form}
\end{equation}
Thus the {\it Calder\'{o}n-Lozanovski{\u \i} space} $E_{\varphi } = \rho _{\varphi }(E, L^{\infty })$ for any Young 
function $\varphi$ is defined by
\begin{equation*}
E_{\varphi } = \{x\in L^{0}:I_{\varphi }(cx)<\infty ~\mathrm{for~some}
~c=c(x)>0\},
\end{equation*}
and it is a Banach ideal space on $\Omega $ with the so-called {\it Luxemburg-Nakano norm}
\begin{equation*}
\| x\|_{E_{\varphi }}=\inf \left\{ \lambda >0:I_{\varphi }\left(
x/\lambda \right) \leq 1\right\} ,
\end{equation*}
where the convex semimodular $I_{\varphi }$ is defined as 
\begin{equation*}
I_{\varphi }(x): = \left\{ 
\begin{array}{cc}
\| \varphi \left( |x| \right) \|_{E} & \text{if }\varphi \left( |x| \right) \in E,\text{ }
\\ 
\infty & \text{otherwise.}
\end{array}
\right.
\end{equation*}
If $E=L^{1}$ ($E=l^{1}$), then $E_{\varphi }$ is the classical {\it Orlicz function}
({\it sequence}) {\it space} $L^{\varphi }$ ($l^{\varphi }$) equipped with the
Luxemburg-Nakano norm (cf. [KR61], [Ma89]). If $E$ is a Lorentz
function (sequence) space $\Lambda _{w}$ ($\lambda _{w}$), then $E_{\varphi
} $ is the corresponding {\it Orlicz-Lorentz function} ({\it sequence}) {\it space} $\Lambda
_{\varphi ,w}$ ($\lambda _{\varphi ,w}$), equipped with the Luxemburg-Nakano
norm. On the other hand, if $\varphi (u)=u^{p},1\leq p<\infty $, then $E_{\varphi }$ 
is the {\it $p$-convexification} $E^{(p)}$ of $E$ with the norm 
$\|x\|_{E^{(p)}}=\| |x|^{p} \|_{E}^{1/p}$. In case $0<p<1,$ we will say
about {\it $p$-concavification} of $E$.

For two ideal (quasi-) Banach spaces $E$ and $F$ on $\Omega $ the symbol 
$E\overset{C}{\hookrightarrow }F$ means that the embedding $E\subset F$ 
is continuous with the norm which is not bigger than $C$, i.e., 
$\| x\|_{F}\leq C\|x\|_{E}$ for all $x\in E$. In the case when the embedding 
$E\overset{C}{\hookrightarrow }F$ holds with some (unknown) constant $C>0$ 
we simply write $E\hookrightarrow F$. Moreover, $E=F$ (and $E\equiv F$) 
means that the spaces are the same and the norms are equivalent (equal).

We will also need some facts from the theory of symmetric spaces. By a 
{\it symmetric function space} (symmetric Banach function space or
rearrangement invariant Banach function space) on $I$, where $I = (0,1)$ or 
$I = (0,\infty )$ with the Lebesgue measure $m$, we mean a Banach ideal space 
$E=(E,\| \cdot \|_{E})$ with the additional property that for any two
equimeasurable functions $x\sim y, x, y\in L^{0}(I)$ (that is, they have the
same distribution functions $d_{x}\equiv d_{y}$, where $d_{x}(\lambda
)=m(\{t\in I:|x(t)|>\lambda \}),\lambda \geq 0$) and $x\in E$ we have that 
$y\in E$ and $\| x\|_{E} = \| y\|_{E}$. In particular, $\| x\|_{E}=\| x^{\ast }\|_{E}$, 
where $x^{\ast }(t)=\mathrm{\inf } \{\lambda >0\colon \ d_{x}(\lambda ) < t\},\ t\geq 0$. 
Similarly, if $e$ is a Banach sequence space with the above property, then 
we say about {\it symmetric sequence space.} It is worth to point out that any Banach
ideal space with this property is equivalent to a symmetric space over one
of the above three measure spaces (cf. [LT79]).

The {\it fundamental function} $f_{E}$ of a symmetric function space $E$
on $I$ is defined by the formula $f_{E}(t)=\| \chi _{\lbrack 0,\,t]}\|_{E},t\in I$. 
It is well-known that each fundamental function is quasi-concave on $I$, that 
is, $f_{E}(0)=0,f_{E}(t)$ is positive, non-decreasing and $f_{E}(t)/t$ is 
non-increasing for $t\in (0,m(I))$ or, equivalently, $f_{E}(t)\leq \max (1,t/s)f_{E}(s)$ 
for all $s,t\in (0, m(I))$.
Moreover, for each fundamental function $f_{E}$, there is an equivalent,
concave function ${\ \tilde{f}_{E},}$ defined by ${\ \tilde{f}_{E}}
(t):=\inf_{s\in (0,m(I))}(1+\frac{t}{s})f_{E}(s)$. Then $f_{E}(t)\leq {
\tilde{f}_{E}}(t)\leq 2f_{E}(t)$ for all $t\in I$. For any quasi-concave
function $\phi $ on $I$ the {\it Marcinkiewicz function space} $M_{\phi }$
is defined by the norm 
\begin{equation*}
\| x \|_{M_{\phi }}=\sup_{ t\in I}\phi (t)\,x^{\ast \ast }(t),~~x^{\ast
\ast }(t)=\frac{1}{t}\int_{0}^{t}x^{\ast }(s)ds.  \label{def marcin}
\end{equation*}
This is a symmetric Banach function space on $I$ with the fundamental
function $f_{M_{\phi }}(t)=\phi (t)$ and $E\overset{1}{\hookrightarrow }
M_{f_{E}}$ since 
\begin{equation}
x^{\ast \ast }(t)\leq \frac{1}{t}\,\| x^{\ast }\|_{E}\| \chi
_{\lbrack 0,t]}\|_{E^{\prime }}=\| x \|_{E}\frac{1}{f_{E}(t)}~
\mathrm{for~any} \,~t\in I,  \label{inequality4}
\end{equation}
(see, for example, [KPS82] or [BS88]). Although the fundamental
function of a symmetric function space $E$ need not be concave, there always
exists equivalent norm on $E$ for which new fundamental function is concave
(cf. Zippin [Zi71], Lemma 2.1). Then for a symmetric function space $E$
with the concave fundamental function $f_{E}$ there is also the smallest
symmetric space with the same fundamental function. This space is the 
{\it Lorentz function space} given by the norm 
\begin{equation*}
\| x \|_{\Lambda _{f_{E}}}=\int_{I}x^{\ast
}(t)df_{E}(t)=f_{E}(0^{+})\| x\|_{L^{\infty }(I)}+\int_{I}x^{\ast
}(t)f_{E}^{\prime }(t)dt.
\end{equation*}
Then the embeddings 
\begin{equation}
\Lambda _{f_{E}}\overset{1}{\hookrightarrow }E\overset{1}{\hookrightarrow }
M_{f_{E}}  \label{estr wloz}
\end{equation}
are satisfied, where $f_{E}$ is the fundamental functions of $E$.

More information about Banach ideal spaces, quasi-Banach ideal spaces, symmetric 
Banach and quasi-Banach spaces can be found, for example, in [KA77], [LT79], [KPR84], 
[JMST], [KPS82] and [BS88].

The paper is organized as follows: In Section 1 some necessary definitions
and notations are collected including the Calder\'{o}n-Lozanovski{\u \i} 
$E_{\varphi }$-spaces. 
In Section 2 the product space $E\odot F$ is defined and some general results 
are presented. We prove important representation of $E\odot F$ as 
$\frac{1}{2}$-concavification of the Calder\'{o}n space $E^{1/2} F^{1/2}$, i.e., 
$E\odot F\equiv \left(E^{1/2} F^{1/2}\right)^{(1/2)}$. Such an equality was used 
by Schep in [Sc10] but without any explanation, which seems to be not 
so evident. Then we present some properties of $E\odot F$ that follow from 
this representation. In particular, the symmetry is proved and formula for the fundamental
function of the product space is given $f _{E\odot F} (t) = f _{E}( t) f _{F}( t) $. 
We finish Section 2 with some sufficient conditions on $E$ and $F$ that $E\odot F$ is 
a Banach space (not only a quasi-Banach space).

In Section 3 we collected properties connecting product spaces with the space 
of multipliers. There is a proof of cancellation property of product operation for 
multipliers $M(E \odot F, E \odot G) \equiv M(F, G)$.

Section 4 is devoted to products of the Calder\'{o}n-Lozanovski{\u \i} spaces of the type 
$E_{\varphi }$ as improvement of results on products known for Orlicz spaces 
$L^{\varphi}$ proved by Ando [An60], Wang [Wa63] and O'Neil [ON65]. 
The inclusion $E_{\varphi _{1}}\odot E_{\varphi_{2}}\hookrightarrow E_{\varphi }$ 
follows from the results proved in [KLM12]. The reverse inclusion $E_{\varphi
}\hookrightarrow E_{\varphi _{1}}\odot E_{\varphi _{2}}$ is investigated here and 
we improve the sufficient and necessary conditions which were given in the case of Orlicz 
spaces by Zabre{\u \i}ko-Ruticki{\u \i} [ZR], Dankert [Da] and Maligranda [Ma89]. 
Combinig the above two inclusions we obtain conditions on equality 
$E_{\varphi } = E_{\varphi_1}\odot E_{\varphi _2}$. For example, for two Young 
functions $\varphi _1,\varphi _2$ we always have 
$E_{\varphi _1} \odot E_{\varphi _2} = E_{\varphi }$, where 
$\varphi =\varphi _1 \oplus \varphi _2$ is defined by 
\begin{equation*}
( \varphi _1 \oplus \varphi _2)( u) = \inf_{u=vw}\, [ \varphi _1( v) + \varphi _2( w)].
\end{equation*}

In Section 5 we deal with the product space of Lorentz and Marcinkiewicz spaces. 
The products of those spaces are calculated. One of the main tools in the proof is 
commutativity of Calder\'on construction with the symmetrizations (cf. Lemma 4).

Section 6 starts with some general discussion about factorization. We prove
that so-called $E$-perfectness of $F$ is necessary for the factorization $F \equiv E\odot M( E, F)$. 
The rest of this section is divided into two parts. The first is devoted to factorization of 
the Calder\'on-Lozanovski{\u \i} $E_{\varphi }$-spaces. Using the results from Section 4 
and the paper [KLM12] we examine when $E_{\varphi }$ can be factorized through 
$E_{\varphi _1}$. In the second part we investigate possibility of factorization for 
Lorentz and Marcinkiewicz spaces. Finally, in Theorem 11 there is proved that under some 
natural assumptions a rearrangement invariant Banach function $X$ space may be factorized 
through Marcinkiewicz space and by duality, Lorentz space may be factorized through a 
rearrangement invariant Banach function space $X$.   

\vspace{3mm}
\begin{center}
\textbf{2. On the product space $E\odot F$}
\end{center}

Given two Banach ideal spaces (real or complex) $E$ and $F$ on 
$( \Omega, \Sigma, \mu)$ define the {\it pointwise product space} $E\odot F$ as 
\begin{equation*}
E\odot F = \left\{ x\cdot y: x\in E ~{\rm and } ~ y\in F\right\} .
\end{equation*}
with a functional $\| \cdot \|_{E\odot F}$ defined by the formula
\begin{equation}
\| z \|_{E\odot F} = \inf \left\{ \| x \|_{E}\, \| y\|_{F}: z = xy, x \in E, y \in F \right \}.  \label{1}
\end{equation}
We will show in the sequel that $E\odot F$ is, in general, a quasi-Banach ideal
space even if both $E$ and $F$ are Banach ideal spaces. 
Let us collect some general properties of the product space and its norm.

\vspace{3mm} 
{\bf Proposition 1.} {\it If $E$ and $F$ are Banach ideal spaces on $( \Omega ,\Sigma ,\mu)$, 
then $E\odot F$ has an ideal property. Moreover, 
\begin{eqnarray*} 
\| z \|_{E\odot F} 
&=& 
 \| \,|z|\, \|_{E\odot F} \\
 &=&
\inf \left\{ \| x \|_{E}\, \| y\|_{F}: |z| = xy, \,x \in E_+, y \in F_+ \right \}\\
&=&
\inf \left\{ \| x \|_{E}\, \| y\|_{F}: |z| \leq xy, \, x \in E_+, y \in F_+ \right \}.
\end{eqnarray*}

}

\begin{proof} We show first that $\| z \|_{E\odot F} = \| \,|z|\, \|_{E\odot F} $. If $z = xy$ with $x \in E, y \in F$, 
then $|z| = z e^{i \theta} = xy e^{i\theta}$, where $\theta: \Omega \rightarrow \mathbb R$, and
$$
\| \,|z|\, \|_{E\odot F} \leq \| x e^{i\theta}\|_E \|y\|_F = \|x\|_E \|y\|_F.
$$
Hence, $\| \,|z|\, \|_{E\odot F} \leq \| z \|_{E\odot F}.$ Similarly, if $|z| = xy$ with $x \in E, y \in F$, then
$z = |z| e^{-i\theta} = xy e^{-i \theta}$ and 
$$
\| z \|_{E\odot F} \leq \| x e^{-i \theta}\|_E \|y\|_F = \| x\|_E \|y\|_F,
$$
from which we obtain the estimate $\| z \|_{E\odot F}  \leq \| \,|z|\, \|_{E\odot F}$. Combining these 
above estimates we obtain $\| z \|_{E\odot F} = \| \,|z|\, \|_{E\odot F} $.

To show the ideal property of $E\odot F$ assume that $z \in E\odot F$ and $|w| \leq |z|$. By definition for any 
$\varepsilon > 0$ we can find $x \in E, y \in F$ such that $z = xy$ and 
$\|x\|_E \|y\|_F \leq \| z \|_{E\odot F} + \varepsilon$. We set $h(t) = \frac{w(t)}{z(t)}$ if $z(t) \neq 0$ and 
$h(t) = 0$ if $z(t) = 0$. Then $w = h z = h x y$ and since $|hx| \leq |x|$ we have $w = hxy \in E\odot F$ with
$$
\| w \|_{E\odot F} \leq \| hx\|_E \|y\|_F \leq \|x\|_E \|y\|_F \leq \| z \|_{E\odot F} + \varepsilon.
$$
Since $ \varepsilon > 0$ was arbitrary, we have $\| w \|_{E\odot F} \leq \| z \|_{E\odot F}$.
Note that in the above proofs we needed only ideal property of one of the spaces $E$ or $F$. 
\vspace{2mm}

Next, if $|z(t)| = x(t)\, y(t), t \in \Omega$, then taking $x_0 = |x|, y_0 = |y|$ we obtain 
$x_0 \geq 0, y_0 \geq 0, x_0\, y_0 = |x|\, |y| = |x y| = |z|$, which gives the proof of the second equality. 
The proof of the third equality follows from the fact that if $0 \leq x \in E, 0 \leq y \in F$ and $|z| \leq xy$, 
then $|z| = u\, x y = x_0 \, y_0$, where $x_0 = u x, y_0 = y$ and $u = \frac{|z|}{xy}$ on the support of 
$x y$ and $u = 0$ elsewhere. 
Since $0 \leq u \leq 1$ it follows that $x_0 \leq x, y_0 \leq y$ and this proves (the non-trivial 
part of) the last equality. 
\end{proof}

Proposition 1 shows that in investigation of product space it is enough to consider real spaces, therefore 
from now we will consider only real Banach ideal spaces.

The product space can be described with the help of the Calder\'{o}n construction. To come to 
this result we first prove some description of $E^{1/p} F^{1-1/p}$ spaces and $p$-convexification. 
Let us start in Theorem 1(i) below with a reformulation of Lemma 31 from [KL10].
\vspace{3mm}
 
{\bf Theorem 1}. {\it Let $E$ and $F$ be a couple of Banach ideal spaces on $(\Omega, \Sigma, \mu) $.
\vspace{-1mm} 

\begin{itemize}
\item[$(i)$] If $1 < p< \infty$ and $z \in E^{1/p} F^{1-1/p}$, then
$$ 
\| z\|_{E^{1/p} F^{1-1/p}}= \inf \left\{ \max \left\{
\| x \|_{E}, \| y \|_{F}\right\}: | z| = x^{1/p}y^{1-1/p}, \, x \in E_+, y \in F_+ \right\} \\
$$
$$
= \inf \left\{ \max \left\{ \| x \|_{E}, \| y \|_{F}\right\}: | z| = x^{1/p}y^{1-1/p},
\| x \| _{E} = \| y \| _{F}, \, x \in E_+, y \in F_+ \right\}.
$$
\item[$(ii)$] If $1 < p< \infty$, then
$$
E^{(p)}\odot F^{(p^{\prime})} \equiv E^{1/p} F^{1-1/p}, ~{\it where} ~ 1/p + 1/p^{\prime} = 1.
$$
\item[$(iii)$] For $0 < p < \infty, ~ (E \odot F)^{(p)} \equiv E^{(p)}\odot F^{(p)}.$
\item[$(iv)$] We have
\begin{equation} \label{representation}
E \odot F \equiv (E^{1/2} F^{1/2})^{(1/2)},
\end{equation}
that is,
\begin{equation} \label{equation8}
\| z \|_{E\odot F}
= \inf \left\{ \max \left\{ \| x \|_{E}^{2}, \| y \|_{F}^{2}\right\}: | z | = xy,
\| x \|_{E} = \| y \|_{F}, \, x \in E_+, y \in F_+ \right\} .
\end{equation}

\end{itemize}
}
\begin{proof} (i) Let $z\in E^{1/p} F^{1-1/p}$ and $z=x^{1/p}y^{1-1/p}$, where 
$0 \leq x \in E, 0 \leq y \in F$. Suppose $\frac{\| x \|_E}{\| y \|_F} = a >1$ (for $0 < a < 1$ proof is similar). 
Put 
\begin{equation*}
x_{1} = a^{-(1-1/p)}\,x,\, y_{1} = a^{\frac{1}{p}} \,y.
\end{equation*}
Then $\| x_{1} \|_{E} = \| x\|_E^{1/p} \| y \|_F^{1-1/p} = \| y_{1} \|_{F}$ and 
$z = x_{1}^{1/p}y_{1}^{1-1/p}$. Of course,  
$$
\max \left\{ \|x_1 \|_E, \,\| y_1 \|_F \right\}  \leq \max \left\{ \| x \|_E,\, \| y \|_F \right\},
$$
which ends the proof. 
\vspace{2mm}

(ii) If $z\in E^{(p)}\odot F^{(p^{\prime})}$, then using Proposition 1 and definition of 
$p$-convexification we obtain
\begin{eqnarray*}
\| z \|_{E^{(p)}\odot F^{(p^{\prime})}}
&=&
\inf \left\{ \| g \|_{E^{(p)}} \| h \|_{F^{(p^{\prime})}}: | z | = gh, \, 0 \leq g \in E^{(p)}, 0 \leq h \in
F^{(p^{\prime})} \right\} \\
&=&
\inf \left\{ \| x \|_{E}^{1/p} \| y \|_{F}^{1-1/p}: | z | = x^{1/p} y^{1-1/p}, \, x \in E_+, y \in F_+ \right\},
\end{eqnarray*}
and using Theorem 1(i) to the last expression we get
\begin{equation*}
\inf_{a>0} \left[ \inf \left \{ \| x \|_{E}^{1/p} \| y \|_{F}^{1-1/p}: | z | = x^{1/p} y^{1-1/p}, \, 
\frac{\| x \|_E}{\| y \|_F} = a, \, x \in E_+, y \in F_+ \right\} \right]
\end{equation*}
\begin{equation*}
= \inf_{a>0}\left[ \inf \left\{ a^{1/p} \| u \|_{E}^{1/p} \| y \|_{F}^{1-1/p}: | z | = a^{1/p}u^{1/p} y^{1-1/p}, 
\| u \|_{E} = \| y \|_{F}, \, u \in E_+, y \in F_+ \right\} \right] 
\end{equation*}
\begin{equation*}
= \inf_{a>0} \left[ a^{1/p} \inf \left\{ \| u \|_E: \frac{| z | }{a^{1/p}} = u^{1/p}y^{1-1/p}, 
\| u \|_E = \| y \| _F, \, u \in E_+, y \in F_+ \right\} \right]
\end{equation*}
\begin{equation*}
=\inf_{a>0}\left[ a^{1/p} \| \frac{z}{a^{1/p}} \|_{E^{1/p} F^{1-1/p}}\right] = \| z \| _{E^{1/p} F^{1-1/p}}.
\end{equation*}

(iii) One has 
\begin{eqnarray*}
\|  z \|_{( E\odot F)^{(p)}}
&=&
\| |z|^p \|_{E\odot F}^{1/p}\\
&=&
\inf \left\{ \| x \|_{E}^{1/p}\, \| y \|_{F}^{1/p}: | z |^p = xy, x \in E_+, y \in F_+\right\} \\
&=&
\inf \left\{ \| u^p \|_{E}^{1/p} \, \| v^p \|_{F}^{1/p}: | z |^p = u^p v^p, u \in E^{(p)}_+, v \in F^{(p)}_+ \right\} \\
&=&
\inf \left\{ \| u \|_{E^{(p)}} \, \| v \|_{F^{(p)}}: | z | = uv, \, u \in E^{(p)}_+, v \in F^{(p)}_+ \right\} \\
&=&
\| z \| _{E^{(p)}\odot F^{(p)}}.
\end{eqnarray*}

(iv) The proof is an immediate consequence of Theorem 1(ii) and (iii) since 
\begin{equation*}
E\odot F\equiv \left( \left( E\odot F\right) ^{(2)}\right) ^{(1/2)} \equiv \left(
E^{1/2} F^{1/2}\right)^{(1/2)}.
\end{equation*}
Moreover, 
\begin{eqnarray*}
\| z \|_{E\odot F}
&=&
\| z \|_{\left(E^{1/2} F^{1/2}\right) ^{(1/2)}} =
( \| \sqrt{ | z |}\,  \|_{E^{1/2} F^{1/2}})^{2}\\
&=&
\left[ \inf \left\{ \max \left\{ \| x \|_{E}, \| y \|_{F}\right\}: \sqrt{| z |} = \sqrt{xy}, 
\| x \| _{E} = \| y \|_{F}, \, x \in E_+, y \in F_+ \right\} \right] ^{2}\\
&=&
\inf \left\{ \max \left\{ \| x \|_{E}^{2}, \| y \|_{F}^{2}\right\}: | z | = x y, \| x \|_{E} = \| y \| _{F}, \,
x \in E_+, y \in F_+ \right\},
\end{eqnarray*}
and the proof is complete.
\end{proof}

As a consequence of the representation (\ref{representation}) we obtain 
the following results:

\vspace{3mm} 
{\bf Corollary 1}. {\it Let $E$ and $F$ be a couple of
Banach ideal spaces on $( \Omega, \Sigma, \mu)$. 
\begin{itemize}
\item[$(i)$] Then $E\odot F$ is a quasi-Banach ideal
space and the triangle inequality is satisfied with constant $2$, i.e.,
\begin{equation*}
\| x+y \|_{E\odot F}\leq 2\left( \| x \|_{E\odot F} + \| y \| _{E\odot F}\right).
\end{equation*}
\item[$(ii)$] If both $E$ and $F$ satisfy the Fatou property, then $E\odot F$ 
has the Fatou property. 
\vspace{-2mm}

\item[$(iii)$] The space $E\odot F$ has order continuous norm if and only if
the couple $(E, F)$ is not jointly order discontinuous, i.e., $( E, F) \not\in ( JOD)$.
\end{itemize}
}

Recall that $( E, F) \in \left( JOD\right)$ (see [KL10]) means that there exist elements 
$x \in E \backslash E_a, y\in F \backslash F_a$ and a sequence of measurable sets 
$A_n\searrow \emptyset$ such that for any sequence $( B_n) $ in 
$\Sigma $ with $B_n\subset A_n$ ($n\in \mathbb{N}$) there are a number 
$a>0$ and a subsequence $(n_{k})$ in $\mathbb{N} $ such that either
\begin{equation*}
\| x \chi_{B_{n_k}} \| _{E}\geq a ~ {\rm and} ~~ \| y\chi _{B_{n_k}} \|_{F}\geq a ~~ {\rm for ~all } ~k \in \mathbb N,
\end{equation*}
or
\begin{equation*}
\| x \chi _{A_{n_k}\backslash B_{n_k}} \|_{E}\geq a ~ {\rm and} ~ 
\| y\chi _{A_{n_k}\backslash B_{n_k}} \|_{F}\geq a ~ {\rm for ~ all} ~ k \in \mathbb N.
\end{equation*}

\begin{proof} (i) It is a consequence of the representation (\ref{representation}) $E\odot
F\equiv ( E^{1/2} F^{1/2})^{(1/2)}$ and the fact that for $1/2$-concavification of a Banach ideal 
space $G = E^{1/2} F^{1/2}$ we have
\begin{eqnarray*}
\| \, |x + y|^{1/2} \, \|_G^2
&\leq&
\left( \| \, |x|^{1/2} \, \|_G + \| \, |y|^{1/2} \, \|_G \right)^2\\
&\leq&
2 \left( \| \, |x|^{1/2} \, \|_G^2 + \| \, | y|^{1/2} \, \|_G^2 \right).
\end{eqnarray*}

(ii) It is again a consequence of the representation (\ref{representation}) and the fact that 
$E^{1/2} F^{1/2}$ has the Fatou property when $E$ and $F$ have the Fatou property 
(see [Lo69], p. 595).
\vspace{2mm}

(iii) The representation (\ref{representation}) and Theorem 13 in [KL10] showing
that $E^{1/2} F^{1/2}$ has order continuous norm if and only if $( E, F) \not\in
\left( JOD\right)$ which gives the statement. 
\end{proof}

Lozanovski{\u \i} [Lo65, Theorem 4] formulated result on the K\"othe dual of $p$-convexification 
$E^{(p)}$ with no proof. The proof can be found in paper by Schep [Sc10, Theorem 2.9] and we 
present another proof which follows from the Lozanovski{\u \i} duality result and our 
Theorem 1(ii). 

\vspace{3mm} 
{\bf Corollary 2}. {\it Let $E$ be a Banach ideal space and $1 < p < \infty$. Then
\begin{equation*}
[E^{(p)}]^{\prime} \equiv (E^{\prime})^{(p)} \odot L^{p^{\prime}}.
\end{equation*}

}

\begin{proof} Using Lozanovski{\u \i} theorem on duality of the Calder\'on spaces 
(see [Lo69], Theorem 2) and our Theorem 1(ii) we obtain
\begin{eqnarray*}
[E^{(p)}]^{\prime} 
&\equiv&
[E^{1/p}(L^{\infty})^{1-1/p}]^{\prime} \equiv (E^{\prime})^{1/p} (L^1)^{1-1/p} \\
&\equiv& 
(E^{\prime})^{(p)} \odot (L^1)^{(p^{\prime})} \equiv (E^{\prime})^{(p)} \odot L^{p^{\prime}} .
\end{eqnarray*}

\vspace{-7mm}

\end{proof}

{\bf Remark 1.} In general, $[E^{(p)}]^{\prime} \neq (E^{\prime})^{(p)}$. In fact, for the classical Lorentz 
space $E = L^{r, 1}$ with $1 < r < \infty$ we have 
$[(L^{r, 1})^{(p)}]^{\prime} = (L^{rp, p})^{\prime} = L^{q, p^{\prime}}$, where $1/q + 1/(pr) =1$ and 
$[(L^{r, 1})^{\prime}]^{(p)} = (L^{r^{\prime}, \infty})^{(p)} = L^{r^{\prime} p, \infty}$.

\vspace{3mm} 
\textbf{Example 1}. (a) If $1 \leq p, q \leq \infty, 1/p + 1/q = 1/r$, then 
$L^p\odot L^q \equiv L^r$. In particular, $L^p\odot L^p \equiv L^{p/2}$.
In fact, by the H\"older-Rogers inequality $\| xy \|_r \leq \| x\|_p \| y\|_q$ for $x \in L^p, y \in L^q$ 
(see [Ma89], p. 69) we obtain $L^p\odot L^q \subset L^r$ 
and $\|z \|_r \leq \|z\|_{L^p\odot L^q}$. On the other hand, if $z \in L^r$, then 
$x = |z|^{r/p} {\rm sgn} z \in L^p, y = |z|^{r/q} {\rm sgn} z \in L^q$ and $x y = z$ with 
$$
\| x\|_p \|y\|_q = \| z\|_r^{r/p} \|z\|_r^{r/q} = \| z\|_r^{r/p+r/q} = \|z\|_r,
$$ 
which shows that $L^r \subset L^p\odot L^q$ and
$ \|z\|_{L^p\odot L^q} \leq \|z\|_r$.

More general, if $1 \leq p, q < \infty, 1/p + 1/q = 1/r$ and $E$ is a Banach ideal space, then 
$E^{(p)} \odot E^{(q)} \equiv E^{(r)}$ (cf. [MP89, Lemma 1] and [ORS08, Lemma 2.21(i)]).
\vspace{2mm}

(b) We have $c_0\odot l^1 \equiv l^{\infty} \odot l^1\equiv l^1$ and 
$c_0 \equiv c_0\odot l^{\infty} \neq l^{\infty}\odot l^{\infty} \equiv l^{\infty}$.
\vspace{3mm}

This example shows that for the Fatou property of $E\odot F$ it is not necessary 
that both $E$ and $F$ do have the Fatou property.

\vspace{3mm}

The next interesting question about product space is its symmetry.
\vspace{3mm}

{\bf Theorem 2}. {\it Let $E$ and $F$ be symmetric Banach spaces on $I = (0, 1)$ or $I = (0, \infty)$ with 
the fundamental functions $f_{E}$ and $f_{F}$, respectively. Then $E\odot F$ is a symmetric quasi-Banach 
space on $I$ and its fundamental function $f_{E\odot F}$ is given by the formula 
\begin{equation} \label{fundamental}
f_{E\odot F}( t) = f_{E}( t) f_{F}( t) ~~ {\rm for} ~ t \in I .
\end{equation}

}

\begin{proof} Using Lemma 4.3 from [KPS82, p. 93] we can easily show that $E^{1/2} F^{1/2}$ 
(even $\rho(E, F)$) is a Banach symmetric space and the representation (\ref{representation}) gives the 
symmetry property of $E\odot F$. 

The inequality $f_{E\odot F}(t) \leq f_{E} (t) f_{F}(t)$ for $t \in I$ is clear. We prove the reverse 
inequality using the fact that each symmetric Banach space $E$ satisfies 
$E\overset{1}{\hookrightarrow }M_{f_{E}},$ where $M_{f_{E}}$ is
the Marcinkiewicz space (see estimate (\ref{inequality4})), some classical inequality on rearrangement 
(see, for example, [HLP52], p. 277 or [BS88], p. 44 or [KPS82], p. 64) and the reverse Chebyshev 
inequality (see Lemma 1 below). For any $0 \leq x \in E, 0 \leq y \in F$ such that $ x\,y = \chi _{[0, t] }$ 
we have
\begin{eqnarray*}
\| x \|_E \, \| y \|_F
&\geq& 
\| x \|_{M_{f_E}} \, \| y \|_{M_{f_F}}\\
&\geq&
\sup_{0 < u \leq t} \frac{f_E ( u)}{u} \int_0^u x^{\ast }( s) \,ds \, \sup_{0 < v \leq t} \frac{f_F ( v)}{v} \int_0^v y^{\ast }( s) \,ds \\ 
&\geq&
\frac{f_E ( t) f_F(t)}{t} \, \frac{1}{t} \int_0^t x^{\ast }(s)\, ds\,  \int_0^t y^{\ast }( s) \,ds \\
&\geq&
\frac{f_E ( t) f_F(t) }{t} \frac{1}{t} \int_0^t x( s) \,ds \, \int_0^t y(s)\, ds\\
&\geq&
\frac{f_E ( t) f_F(t) }{t} \int_0^t x(s) \, y(s)\, ds\\
&=&
 \frac{f_E ( t) f_F(t) }{t} \int_0^t \chi_{[0, t]}(s)\, ds = f_{E} ( t) f_F(t),
\end{eqnarray*}
and so $\| \chi _{[0, t]} \| _{E\odot F} \geq f_{E} ( t) f_F(t)$.
\end{proof}

In the fifth inequality above we used the reverse Chebyshev inequality which we will prove in the lemma 
below. On the classical Chebyshev inequality for decreasing functions (see, for example, [Mi70], p. 39 or [HM91-2], p. 213).

\vspace{3mm}
\textbf{Lemma 1.} \label{Lemma1} {\it Let $0 < \mu(A) < \infty$ and $x(s)\, y(s) = a > 0$ for all $s \in A$
with $0 \leq x, y \in L^1(A)$. Then
\begin{equation} \label{equation9}
\mu(A) \int_A x y\, d\mu \leq \int_A x\,d\mu\, \int_A y\,d\mu.
\end{equation} 

}

\begin{proof} For any $s, t \in A$ we have
\begin{eqnarray*}
x(s) y(s) &-& x(s) y(t) - x(t) y(s) + x(t) y(t) 
=
2 a - x(s) y(t) - x(t) y(s)\\
&=& 
2 a - \frac{x(s) a}{x(t)} - \frac{x(t) a}{x(s)} = a \, [2 - \frac{x(s)}{x(t)} - \frac{x(t)}{x(s)}]\\
&=&
a\, \frac{2x(s)x(t) - x(s)^2 - x(t)^2}{x(s)x(t)} = - a \, \frac{[x(s) - x(t)]^2}{x(s)x(t)} \leq 0.
\end{eqnarray*}
Now integrating over $A$ with respect to $s$ and over $A$ with respect to $t$ we obtain 
the desired inequality (\ref{equation9}).
\end{proof}

{\bf Remark 2.} Formula (\ref{fundamental}) is a generalization of the well-known equality 
on fundamental functions $f_{E}( t) f_{E^{\prime }}( t) = t = f_{L^{1}}( t) $ for $t\in I$ and it 
is also true for symmetric sequence spaces with the same proof. 
\vspace{2mm}

Example 1(a) shows that $E\odot F$ is, in general, a quasi-Banach ideal
space even if both $E$ and $F$ are Banach ideal spaces. We can ask under which 
additional conditions on $E$ and $F$ the product space $E\odot F$ is a Banach ideal 
space. Before formulation the theorem we need notion of $p$-convexity.
A Banach lattice $E$ is said to be {\it $p$-convex} ($1 \leq p < \infty$) with constant $K\geq 1$ if
\begin{equation*}
\| ( \sum_{k=1}^{n}  | x_k |^p )^{1/p}\|_E \leq K \, ( \sum_{k=1}^{n} \| x_k \|_E^{p})^{1/p}. 
\end{equation*}
for any sequence $\left(x_k \right)_{k=1}^n \subset X$ and any $n \in \mathbb N$. 
If a Banach lattice $E$ is $p$-convex with constant $K \geq 1$ and $1 \leq q < p$, then 
$E$ is also $q$-convex with constant at most $K$.
Moreover, $p$-convexification $E^{(p)}$ of a Banach lattice $E$ is $p$-convex with constant $1$. 
More information on $p$-convexity we can find, for example, in [LT79] and [Ma04].

\vspace{3mm}
{\bf Theorem 3.} {\it Suppose that $E, F$ are Banach ideal spaces such that $E$ is $p_0$-convex 
with constant $1$, $F$ is $p_1$-convex with constant $1$ and $\frac{1}{p_0}+\frac{1}{p_1} \leq 1$. 
Then $E\odot F$ is a Banach space which is even $\frac{p}{2}$-convex, where 
$\frac{1}{p} = \frac{1}{2}(\frac{1}{p_0} + \frac{1}{p_1})$.}
\vspace{3mm}

Before the proof of Theorem 3 let us present the following lemma, which was mentioned in [Re80] and 
proved in [TJ89], p. 219. For the sake of completeness we give its proof.

\vspace{3mm} 
\textbf{Lemma 2.} \label{Lemma2} {\it If $E$ is $p_0$-convex with constant $K_0$ and $F$ is $p_1$-convex 
with constant $K_1$, then $E^{1-\theta} F^{\theta}$ is $p$-convex with constant 
$K \leq K_0^{1-\theta} K_1^{\theta}$, where $\frac{1}{p} = \frac{1-\theta}{p_0} + \frac{\theta}{p_1}$.}

\begin{proof} Let $z_1, \ldots, z_n \in E^{1-\theta} F^{\theta}$. There are $x_k \in E, y_k \in F$ with 
$\| x_k\|_E \leq 1, \|y_k \|_F \leq 1$ and $\lambda_k \geq 0$ be such that for all $1 \leq k \leq n$
$$
|z_k| \leq \lambda_k |x_k|^{1-\theta} |y_k|^{\theta} ~~ {\rm and} ~~ \lambda_k \leq \| z_k\|_{E^{1-\theta} F^{\theta}} (1 + \varepsilon).
$$ 
By the H\"older-Rogers inequality
\begin{eqnarray*}
( \sum_{k=1}^n |z_k|^{p} )^{1/p} 
&\leq& 
( \sum_{k=1}^n \lambda_k^p\,  |x_k|^{p(1-\theta)} |y_k|^{p \theta} )^{1/p} \\
&\leq&
( \sum_{k=1}^n \lambda_k^{p} \,|x_k|^{p_0})^{(1-\theta)/p_0} ( \sum_{k=1}^n \lambda_k^{p}\, |y_k|^{p_1} )^{\theta/p_1}.
\end{eqnarray*}
Then, by the property $\| u^{1-\theta} v^{\theta} \|_{E^{1-\theta} F^{\theta}} \leq \| u\|_E^{1- \theta} \, \| v\|_F^{\theta}$ and assumptions on convexity together with assumptions on elements, we obtain
\begin{eqnarray*}
 \| ( \sum_{k=1}^n |z_k|^{p} )^{1/p} \|_{E^{1-\theta} F^{\theta}}
&\leq& 
 \| ( \sum_{k=1}^n \lambda_k^{p}\, |x_k|^{p_0})^{(1-\theta)/p_0} ( \sum_{k=1}^n \lambda_k^{p}\, |y_k|^{p_1} )^{\theta/p_1} \|_{E^{1-\theta} F^{\theta}}\\
&\leq&
 \| ( \sum_{k=1}^n \lambda_k^{p} \, |x_k|^{p_0} )^{1/p_0} \|_E^{1-\theta} \,  \| ( \sum_{k=1}^n \lambda_k^{p}\, |y_k|^{p_1} )^{1/p_1} \|_F^{\theta}\\
&\leq&
K_0^{1-\theta} (\sum_{k=1}^n \|  \lambda_k^{p/p_0} \, x_k\|_E^{p_0} )^{(1-\theta)/{p_0}} \, K_1^{\theta} (\sum_{k=1}^n \|  \lambda_k^{p/p_1} \, y_k\|_F^{p_1} )^{\theta/{p_1}}\\
&=&
K_0^{1-\theta} K_1^{\theta} \, (\sum_{k=1}^n  \lambda_k^{p}\, \| x_k\|_E^{p_0} )^{(1-\theta)/{p_0}} \, (\sum_{k=1}^n   \lambda_k^{p}\, \| y_k\|_F^{p_1} )^{\theta/{p_1}}\\
&\leq&
K_0^{1-\theta} K_1^{\theta} \, (\sum_{k=1}^n  \lambda_k^{p} )^{(1-\theta)/{p_0}} \, (\sum_{k=1}^n \lambda_k^{p}  )^{\theta/{p_1}} = K_0^{1-\theta} K_1^{\theta} \, (\sum_{k=1}^n  \lambda_k^{p} )^{1/{p}}\\
&\leq&
(1+\varepsilon) \, K_0^{1-\theta} K_1^{\theta} \,  (\sum_{k=1}^n \| z_k\|_{E^{1-\theta} F^{\theta}}^{p} )^{1/{p}}. 
\end{eqnarray*}
Since $\varepsilon > 0$ is arbitrary, the $p$-convexity of $E^{1-\theta} F^{\theta}$ is proved with the constant $K \leq  K_0^{1-\theta} K_1^{\theta}$.
\end{proof}

\noindent
{\it Proof of Theorem 3}. By Lemma 2, we have that $Z = E^{1/2} F^{1/2}$ is $p$-convex with constant $1$, where
$\frac{1}{p} = \frac{1}{2 p_0} +  \frac{1}{2 p_1}$. The assumption on $p_0, p_1$ gives that $p \geq 2$ and since 
$1/2$-concavification of $Z$ is $p/2$-convex with constant $1$ ($p/2 \geq 1$) it follows that it is $1$-convex with 
constant $1$ which gives that the norm of $E\odot F = Z^{(1/2)}$ satisfies the triangle inequality, and consequently 
is a Banach space. This completes the proof. \qed
\vspace{2mm}

{\bf Remark 3.} By duality arguments, Theorem 3 can be also formulated in the terms of $q$-concavity of 
the K\"othe dual spaces. A Banach lattice $E$ is {\it $q$-concave} ($1 < q < \infty$) with constant $K \geq 1$ 
if $( \sum_{k=1}^n \| x_k\|_F^q)^{1/q} \leq K \, \| (\sum_{k=1}^n |x_k|^q)^{1/q}\|_F$ for any sequence 
$\left(x_k \right)_{k=1}^n \subset X$ and any $n \in \mathbb N$. 
\vspace{2mm}

{\bf Remark 4.} Since imbedding $G \subset E\odot F$ means also factorization $ z = x\,y$, where $x \in E$ and $y \in F$, therefore sometimes these imbeddings or identifications of product spaces $E\odot F = G$ are called factorizations of concrete spaces as, for example, $l^p$ and Ces\`aro sequence spaces or $L^p$ and Ces\`aro function spaces (cf. [Be96], [AM09], [Sc10]), factorization of tent spaces or other spaces (cf. [CV00], [CRW76], [Ho77]).

\vspace{3mm}
\begin{center}
\textbf{3. The product spaces and multipliers}
\end{center}

Let us collect properties connecting product space with the space of multipliers. 
We start with the Cwikel and Nilsson result [CN03, Theorem 3.5]. They proved that if a Banach ideal 
space $E$ has the Fatou property and $0 < \theta < 1$, then
\begin{equation*} \label{CN03}
E \equiv M( F^{(1/\theta)}, E^{1-\theta} F^{\theta})^{(1-\theta)}.
\end{equation*}
We will prove a generalization of this equality, which in the case of $G = L^{\infty}$ coincides with their result.

\vspace{3mm}
{\bf Proposition 2.} {\it Let $E, F, G$ be Banach ideal spaces. Suppose that $E$ has the Fatou property 
and $0 < \theta < 1$. Then 
\begin{equation*} \label{KLM}
M(G, E) \equiv M( G^{1-\theta} F^{\theta}, E^{1-\theta} F^{\theta})^{(1-\theta)}.
\end{equation*}

}

\begin{proof} First, let us prove the imbedding $\overset{1}\hookrightarrow$. Let $x \in M(G, E)$. We want to show that 
$x \in M( G^{1-\theta} F^{\theta}, E^{1-\theta} F^{\theta})^{(1-\theta)}$, that is, $|x|^{1-\theta} \in 
M( G^{1-\theta} F^{\theta}, E^{1-\theta} F^{\theta})$ or equivalently $x^{1-\theta} |y| \in  E^{1-\theta} F^{\theta}$ 
for any $y \in G^{1-\theta} F^{\theta}$. Take arbitrary $y \in G^{1-\theta} F^{\theta}$ with the norm $< 1$. 
Then there are $w \in G, v \in F$ satisfying $\| w\|_G \leq 1, \| v\|_F \leq 1$ and $|y| \leq |w|^{1-\theta} |v|^{\theta}$. 
Clearly,
$$
|x|^{1-\theta} \, |y| \leq | x w|^{1-\theta} |v|^{\theta} \in E^{1-\theta} F^{\theta}
$$
since $x \in M(G, E), w \in G$ gives $x w \in E$. This proves the inclusion. Moreover,
\begin{eqnarray*}
\| |x|^{1-\theta} \, |y| \|_{E^{1-\theta} F^{\theta}}
&\leq&
\|  | x w|^{1-\theta} |v|^{\theta} \|_{E^{1-\theta} F^{\theta}} \\
&\leq&
\|  x w \|_E^{1-\theta} \| v \|_F^{\theta} \\
&\leq&
\|  x \|_{M(G, E)}^{1-\theta} \| w \|_G^{1-\theta} \| v \|_F^{\theta} \leq \| x \|_{M(G, E)}^{1-\theta}.
\end{eqnarray*}
Thus,
$$
\| |x|^{1-\theta} \|_{M(G^{1-\theta} F^{\theta}, \, E^{1-\theta} F^{\theta})}^{1/(1-\theta)} \leq \| x \|_{M(G, E)}.
$$
The imbedding $\overset{1}\hookleftarrow$. Let $\| x \|_{ M( G^{1-\theta} F^{\theta}, E^{1-\theta} F^{\theta})^{(1-\theta)}} = 1$, 
i.e. $\| |x|^{1-\theta}\|_{ M( G^{1-\theta} F^{\theta}, E^{1-\theta} F^{\theta})} = 1$. We need to show that for any 
$w \in G$ we have $x w \in M( F^{(1/\theta)}, E^{1-\theta} F^{\theta})^{(1-\theta)}$, that is, 
$|x w|^{1-\theta} \in M( F^{(1/\theta)}, E^{1-\theta} F^{\theta})$. Really, by the Cwikel-Nilsson result, we 
obtain $x w \in E$ for any $w \in G$.

Let $w \in G$ and $v \in F$. Since the norm of $x$ is 1 it follows that 
$$
\| |x|^{1-\theta} \frac{z}{\| z\|_{ G^{1-\theta} F^{\theta}}} \|_{ E^{1-\theta} F^{\theta}} \leq 1
$$
for each $0 \neq z \in  G^{1-\theta} F^{\theta}$. Consequently, for $ z = |w|^{1-\theta} |v|^{\theta}$, we obtain
\begin{eqnarray*}
\| |x w|^{1-\theta} \, |v|^{\theta} \|_{E^{1-\theta} F^{\theta}}
&=&
\|  | x|^{1-\theta} |w|^{1-\theta} |v|^{\theta} \|_{E^{1-\theta} F^{\theta}} \\
&\leq&
\|  |w|^{1-\theta} |v|^{\theta} \|_{G^{1-\theta}F^{\theta}} \leq \|  w \|_G^{1-\theta} \| v \|_F^{\theta}.
\end{eqnarray*}
This proves our inclusion part because from the assumption on $x$ and the fact that 
$ |w|^{1-\theta} \, |v|^{\theta} \in G^{1-\theta} F^{\theta}$ we have
$ |x w|^{1-\theta} \, |v|^{\theta} =  |x|^{1-\theta} |w|^{1-\theta} \, |v|^{\theta} \in E^{1-\theta} F^{\theta}$. Moreover, 
by the Cwikel-Nilsson result and the last estimate with $v = |m|^{1/\theta}$, we obtain
\begin{eqnarray*}
\| x \|_{M(G, E)} 
&=&
\sup_{\| w\|_G \leq 1} \| x w \|_E = \sup_{\| w\|_G \leq 1} \| x w \|_{M( F^{(1/\theta)}, E^{1-\theta} F^{\theta})^{(1-\theta)}} \\
&=&
\sup_{\| w\|_G \leq 1} \sup_{\| m\|_{F^{(1/\theta)}} \leq 1} \| |x w|^{1-\theta} m \|_{E^{1-\theta} F^{\theta}}^{1/(1-\theta)}\\
&\leq&
\sup_{\| w\|_G \leq 1} \sup_{\| m\|_{F^{(1/\theta)}} \leq 1} \| w \|_G \| m\|_{F^{(1/\theta)}}^{1/(1-\theta)} \leq 1,
\end{eqnarray*}
and the theorem is proved with the equality of the norms.
\end{proof}
 
Proposition 2 together with the representation of the product space as the $1/2$-concavification of the 
Calder\'on space will give the ``cancellation" property for multipliers of products.
 
\vspace{2mm} 
{\bf Theorem 4.} \label{Theorem4} {\it Let $E, F, G$ be Banach ideal spaces. If $G$ has the Fatou property, then 
\begin{equation} \label{cancellation}
M(E \odot F, E \odot  G) \equiv M( F, G).
\end{equation}

}

\begin{proof} Applying Theorem 1(iv), property (g) from [MP89] and Proposition 2 we obtain
\begin{eqnarray*}
M(E \odot F, E \odot  G) 
&\equiv & 
M[ (E^{1/2} F^{1/2})^{(1/2)}, (E^{1/2} G^{1/2})^{(1/2)}] \\
&\equiv &
M[E^{1/2} F^{1/2}, (E^{1/2} G^{1/2})]^{(1/2)}  \equiv M(F, G),
\end{eqnarray*}
\vspace{-2mm}
and (\ref{cancellation}) is proved.
\end{proof}

{\bf Remark 5.} Note that Proposition 2 and Theorem 4 are equivalent. Proposition 2 can be written in the following 
form: if $F$ has the Fatou property, then  
\begin{equation*}
M(E^{\theta} G^{1-\theta}, F^{\theta} G^{1-\theta}) \equiv M(E, F)^{(1/\theta)},
\end{equation*}
and it can be proved using Theorem 4. In fact, applying Theorem 1(ii), cancellation property 
from Theorem 4 and property (g) from [MP89] we obtain
\begin{eqnarray*}
M(E^{\theta} G^{1-\theta}, F^{\theta} G^{1-\theta}) 
&\equiv& 
M(E^{(1/\theta)} \odot G^{(1/(1-\theta))}, F^{(1/\theta)} \odot G^{(1/(1-\theta))}) \\
&\equiv&
M(E^{(1/\theta)}, F^{(1/\theta)}) \equiv M(E, F)^{(1/\theta)}.
\end{eqnarray*}

From Theorem 4 we can also get the equality mentioned by Raynaud [Ra92] which can be proved also directly 
(cf. also [Sc10], Proposition 1.4).
 
\vspace{3mm} 
{\bf Corollary 3.} \label{Corollary6} {\it Let $E, F$ be Banach ideal spaces. If $E$ has the Fatou property, then 
\begin{equation} \label{duality}
(E \odot F)^{\prime} \equiv M(F, E^{\prime}) \equiv M( E, F^{\prime})
\end{equation}

}

\begin{proof} Using the Lozanovski{\u \i} factorization theorem (for more discussion see Part 6) and the 
cancellation property (\ref{cancellation}) we obtain
$$
(E \odot F)^{\prime} \equiv M(E \odot F, L^1) \equiv M(E \odot F, E \odot E^{\prime})  \equiv M(F, E^{\prime}).
$$
and
$$
(E \odot F)^{\prime} \equiv M(E \odot F, L^1) \equiv M(E \odot F, F^{\prime} \odot F)  \equiv M(E, F^{\prime}).
$$
Note that the second identity in (\ref{duality}) follows also from the general properties of multipliers 
(see [MP89, property (e)] or [KLM12, property (vii)]) because we have 
$M(F, E^{\prime}) \equiv M(E^{\prime \prime}, F^{\prime}) \equiv M(E, F^{\prime})$. \end{proof} 
 
\vspace{2mm} 
{\bf Corollary 4.} \label{Corollary4} {\it Let $E, F$ be Banach ideal spaces. If $F$ has the Fatou property and 
$\| xy \|_{E \odot F} \leq 1$ for all $x \in E$ with $\| x\|_E \leq 1$, then $\| y \|_F \leq 1$.}

\begin{proof} Since by assumption $\| y \|_{M(E, E \odot F)} \leq 1$, then using Theorem 4, together with the facts 
that $M(L^{\infty}, F) \equiv F, E \odot L^{\infty} \equiv E$, we obtain 
$$
\| y \|_F = \| y \|_{M(L^{\infty}, F)} = \| y \|_{M(E \odot L^{\infty}, E \odot F)} = \| y \|_{M(E, E \odot F)} \leq 1.
$$

\vspace{-5mm} 

\end{proof} 

{\bf Corollary 5.} \label{Corollary5} {\it Let $E, F, G$ be Banach ideal spaces. If $F$ and $G$ have the Fatou 
property, then 
\begin{equation*}
M(E \odot F, G) \equiv M(E, M(F, G)).
\end{equation*}
}

\vspace{-5mm}

\begin{proof} Using Theorem 4, the Lozanovski{\u \i} factorization theorem, Corollary 3 with the fact that the 
Fatou property of $F$ gives by Corollary 1(ii) that $F \odot G^{\prime}$ has the Fatou property, again 
Corollary 3 and the Fatou property of $G$ we obtain
\begin{eqnarray*}
M(E \odot F, G) 
&\equiv& M(E \odot F \odot G^{\prime}, G \odot G^{\prime}) \\
&\equiv&  
M(E \odot F \odot G^{\prime}, L^1) \\
&\equiv&
(E \odot F \odot G^{\prime})^{\prime} \equiv (F \odot G^{\prime} \odot E)^{\prime}\\
&\equiv&
M(E, (F \odot G^{\prime} )^{\prime}) \equiv M(E, M(F, G^{\prime \prime}) \equiv M(E, M(F, G)),
\end{eqnarray*}
which establishes the formula.
\end{proof} 

\begin{center}
\textbf{4. The product of Calder\'{o}n-Lozanovski{\u \i} $E_{\varphi}$-spaces}
\end{center}

The pointwise product of Orlicz spaces was investigated already by Krasnoselski{\u \i} and Ruticki{\u \i} in their 
book, where sufficient conditions on imbedding $L^{\varphi_1} \odot L^{\varphi_2} \subset L^{\varphi}$ are given 
in the case when $\Omega$ is bounded closed subset of ${\mathbb R}^n$ (cf. [KR61], Theorems 13.7 
and 13.8). For the same set $\Omega$,  Ando [An60] proved that $L^{\varphi_1} \odot L^{\varphi_2} \subset L^\varphi $  
if and only if there exist $C>0$, $u_0>0$ such that $\varphi(C u v)\leq \varphi_1(u)+\varphi_2(v)$ for $u, v \geq u_0$. 

O'Neil [ON65] presented necessary and sufficient conditions for the imbedding $L^{\varphi_1} \odot L^{\varphi_2} 
\subset L^{\varphi}$ in the case when measure space is either non-atomic and infinite or non-atomic and finite 
or counting measure on $\mathbb N$. Moreover, he observed that condition $\varphi(C u v)\leq \varphi_1(u)+\varphi_2(v)$ 
for all [large, small] $u, v > 0$ is equivalent to condition on inverse functions  
$C_1 \varphi _{1}^{-1}( u) \varphi _{2}^{-1}( u) \leq \varphi ^{-1}( u) $ for all [large, small] $u > 0$. O'Neil's results were 
also presented, with his proofs, in the books [Ma89, pp. 71-75] and [RR91, pp. 179-184].

The reverse imbedding $L^\varphi \subset L^{\varphi_1} \odot L^{\varphi_2}$ and the equality $L^{\varphi_1} \odot L^{\varphi_2} = L^\varphi $ were considered by Zabreiko-Ruticki{\u \i} [ZR67, Theorem 8], Dankert [Da74, pp. 63-68] and Maligranda [Ma89, 69-71].

We will prove the above results for more general spaces, that is, for the Calder\'{o}n-Lozanovski{\u \i} $E_{\varphi}$-spaces. 
Results on the imbedding $E_{\varphi_1}\odot E_{\varphi _2}\hookrightarrow E_{\varphi }$ need the following relations 
between Young functions (cf. [ON65]): we say $\varphi _{1}^{-1}\varphi _{2}^{-1}\prec \varphi ^{-1}$ for all arguments [for 
large arguments] (for small arguments) if that there is a constant $C>0$ [there are constants 
$C, u_0 > 0$] (there are constants $C, u_0 > 0$) such that the inequality
\begin{equation} \label{left-ineq}
C\varphi _{1}^{-1}( u) \varphi _{2}^{-1}( u) \leq \varphi ^{-1}( u) 
\end{equation}
holds for all $u>0$ [for all $u\geq u_0$] (for all $u\leq u_0$), respectively.
\vspace{2mm}

{\bf Remark 6.} The inequality (\ref{left-ineq}) implies a generalized Young inequality: 
\begin{equation} \label{Young}
\varphi ( Cuv ) \leq \varphi _{1}(u)+\varphi _{2}(v) ~~ {\rm for ~all} ~~ u,v>0 ~~ {\rm such ~ that} 
~ \varphi _{1} (u) ,\varphi _{2}(v) <\infty .
\end{equation}
On the other hand, if $\varphi (Cuv)\leq \varphi _{1}(u)+\varphi _{2}(v)$ for all $u,v>0$, then 
$\varphi _{1}^{-1}(w)\varphi _{2}^{-1}(w)\leq \frac{2}{C}\varphi ^{-1}(w)$ for each $w>0$ 
(see [ON65] and [KLM12]). Similar equivalences hold for large and small arguments.
\vspace{3mm}

In [KLM12] the question when the product $x y \in E_{\varphi}$ provided $x \in
E_{\varphi_1}$ and $y \in E_{\varphi_2}$ was investigated, as a generalization of O'Neil's 
theorems [ON65], and the following results were proved (see [KLM12], Theorems 4.1, 4.2 and 4.5):

\vspace{3mm}
{\bf Theorem A}. {\it Let  $\varphi _1, \varphi_2$ and $\varphi$ be three Young functions.}
\vspace{-2mm}

\begin{itemize}
\item[$(a)$] {\it If $E$ is a Banach ideal space with the Fatou property and one of the following conditions holds: }

$(a1)$ ~ {\it $\varphi _{1}^{-1}\varphi _{2}^{-1}\prec \varphi^{-1}$
for all arguments,}

$(a2)$ ~ {\it $\varphi _{1}^{-1}\varphi _{2}^{-1} \prec \varphi^{-1} 
$ for large arguments and $L^{\infty }\hookrightarrow E$, }

$(a3)$ ~ {\it $\varphi _{1}^{-1}\varphi _{2}^{-1} \prec
\varphi^{-1} $ for small arguments and $E\hookrightarrow L^{\infty }$, }
\end{itemize}
\vspace{-2mm}

{\it then, for every $x \in E_{\varphi_1}$ and $y \in E_{\varphi_2}$ the product 
$x y \in E_{\varphi}$, which means that 

$ E_{\varphi _1} \odot E_{\varphi _2} \hookrightarrow E_{\varphi}.$ }
\vspace{-2mm}

\begin{itemize}
\item[$(b)$] {\it If a Banach ideal space $E$ with the Fatou property is such that $E_{a}\neq \left\{ 0\right\} $ and we have 
$E_{\varphi _1}\odot E_{\varphi_2}\hookrightarrow E_{\varphi }$, then $\varphi _{1}^{-1}\varphi
_{2}^{-1}\prec \varphi ^{-1}$ for large arguments.}
\vspace{-2mm}

\item[$(c)$] {\it If a Banach ideal space $E$ has the Fatou property, ${\rm supp}E_{a}=\Omega, 
L^{\infty }\not\hookrightarrow E$ and 
$E_{\varphi _1}\odot E_{\varphi _2}\hookrightarrow E_{\varphi }$, 
then $\varphi _{1}^{-1}\varphi _{2}^{-1}\prec \varphi ^{-1}$ for small
arguments.}
\vspace{-2mm}

\item[$(d)$] {\it If $e$ is a Banach sequence space with the Fatou property, 
$\sup_{k \in \mathbb N} \| e_k \|_e < \infty, l^{\infty }\not\hookrightarrow e$ and 
$e_{\varphi _{1}}\odot e_{\varphi_{2}}\hookrightarrow e_{\varphi }$, then 
$\varphi _{1}^{-1}\varphi_{2}^{-1}\prec \varphi ^{-1}$ for small arguments.}
\vspace{-2mm}
\end{itemize}

Note that in the case (c) we can even conclude the relation $\varphi _{1}^{-1}\varphi _{2}^{-1}\prec \varphi ^{-1}$ 
for all arguments, using (b) and (c).
\vspace{2mm}

The sufficient and necessary conditions on the reverse inclusion 
$E_{\varphi } \hookrightarrow E_{\varphi _1}\odot E_{\varphi
_2}$ need also the reverse relations between Young functions, the same as in [KLM12].

The symbol $\varphi ^{-1}\prec \varphi _{1}^{-1}\varphi _{2}^{-1}$ for all arguments 
[for large arguments] (for small arguments) means that there is a constant $D>0$ [there 
are constants $D,u_0 >0$] (there are constants $D,u_0 >0$) such that the inequality 
\begin{equation}
\varphi ^{-1}( u) \leq D \, \varphi_{1}^{-1}( u) \varphi_{2}^{-1} ( u)
\label{right-ineq}
\end{equation}
holds for all $u>0$ [for all $u\geq u_0$] (for all $0 < u\leq u_0$),
respectively.

\vspace{3mm} 
{\bf Theorem 5}. {\it Let  $\varphi _1, \varphi_2$ and $\varphi$ be 
three Young functions.}
\vspace{-2mm}

\begin{itemize}
\item[$(a)$] {\it If $E$ is a Banach ideal space with the Fatou property and one of the following conditions holds: }

$(a1)$ ~ {\it $\varphi^{-1} \prec \varphi _{1}^{-1}\varphi _{2}^{-1}$
for all arguments,}

$(a2)$ ~ {\it $\varphi^{-1} \prec \varphi _{1}^{-1}\varphi _{2}^{-1}$ for large arguments 
and $L^{\infty }\hookrightarrow E$,} 

$(a3)$ ~ {\it $\varphi^{-1} \prec \varphi _{1}^{-1}\varphi _{2}^{-1}$ for small arguments 
and $E\hookrightarrow L^{\infty }$,}
\end{itemize}
\vspace{-2mm}

{\it then $E_{\varphi} \hookrightarrow E_{\varphi _1} \odot E_{\varphi _2}$.}
\vspace{-2mm}

\begin{itemize}
\item[$(b)$] {\it If $E$ is a symmetric Banach function space on $I$ with the Fatou property, $E_{a}\neq \{ 0\} $ and 
$E_{\varphi }\hookrightarrow E_{\varphi _1}\odot E_{\varphi _2}$, 
then $\varphi^{-1} \prec \varphi _{1}^{-1}\varphi _{2}^{-1}$ for large arguments.}
\vspace{-2mm}

\item[$(c)$] {\it If $E$ is a symmetric Banach function space on $I$ with the Fatou property, 
${\it supp} E_{a}=\Omega, L^{\infty}\not\hookrightarrow E$ and $E_{\varphi }\hookrightarrow E_{\varphi
_1}\odot E_{\varphi _2}$, then $\varphi^{-1} \prec \varphi _{1}^{-1}\varphi _{2}^{-1}$ for 
small arguments.}
\vspace{-2mm}

\item[$(d)$] {\it Let $e$ be a symmetric Banach sequence space with the Fatou property and 
order continuous norm. If $e_{\varphi }\hookrightarrow e_{\varphi _{1}}\odot e_{\varphi _{2}}$, 
then $\varphi^{-1} \prec \varphi _{1}^{-1}\varphi _{2}^{-1}$ for small arguments.}
\end{itemize}

\begin{proof} $(a1)$ The idea of the proof is taken from [Ma89], Theorem 10.1(b). For 
$z\in E_{\varphi }\backslash \left\{ 0\right\}$ let $y=\varphi ( \frac{| z |}{\| z \|_{E_{\varphi }}})$ and
\begin{equation*}
z_{i}( t)
=\left\{ 
\begin{array}{ccc}
\sqrt{\frac{| z(t)| }{\varphi _{1}^{-1}( y(t))\, \varphi _{2}^{-1}( y(t)) }}\, \varphi _{i}^{-1}( y( t)),  & \text{if} & 
t\in ${\it supp} z$, \\ 
0, &  & \text{otherwise,}
\end{array}
\right.
\end{equation*}
for $i=1, 2$. The elements $z_{i}$ are well defined. Indeed, if $a_{\varphi }=0$, then 
$y( t) >0$ for $\mu$-a.e. $t\in {\it supp} z$. 
If $a_{\varphi }>0$, then assumption on functions implies that $a_{\varphi _{1}}>0$ 
and $a_{\varphi _{2}}>0$. Consequently, $\varphi_{1}^{-1}( 0) = a_{\varphi _{1}}$ 
and $\varphi _{2}^{-1}(0) =a_{\varphi _{2}}$. Now we will prove the inequality
\begin{equation}\label{equation12}
\varphi _{i} ( \frac{z_i}{\sqrt{D \| z \|_{E_{\varphi }}}}) \leq y, ~  i=1,2.  
\end{equation}
If $a_{\varphi }>0$, taking $u\rightarrow 0$ in inequality (\ref {right-ineq}) we obtain
$a_{\varphi }\leq Da_{\varphi _{1}}a_{\varphi _{2}}.$ If $y(t) =0$, then
\begin{equation*}
z_i(t) = \sqrt{\frac{| z( t) |}{a_{\varphi _1} \, a_{\varphi _2}}} \, \varphi _{i}^{-1}( 0) \leq 
\sqrt{\frac{\| z \|_{E_{\varphi }} a_{\varphi }} {a_{\varphi_{1}} a_{\varphi _{2}}}} \, \varphi _{i}^{-1}( 0) 
\leq \sqrt{D \|z \|_{E_{\varphi }}} \, \varphi _{i}^{-1}( 0)
\end{equation*}
and consequently $\varphi _i ( \frac{z_i( t) }{\sqrt{D \| z \|_{E_{\varphi }}}}) = 0 = y( t).$
If $y\left( t\right) >0$, then
\begin{equation*}
z_i( t) = \sqrt{\frac{ | z ( t)| }{\varphi _{1}^{-1} ( y( t)) \varphi _{2}^{-1} (y( t) ) }} \, \varphi _{i}^{-1}( y( t))
\leq \sqrt{\frac{D | z(t)| }{\varphi^{-1} ( y( t)) }} \, \varphi _{i}^{-1}( y(t) ) =
\sqrt{D \| z \| _{\varphi }} \, \varphi_{i}^{-1}( y(t)).
\end{equation*}
This proves (\ref{equation12}) and consequently we obtain
\begin{equation*}
I_{\varphi _1} ( \frac{z_{1}}{\sqrt{D \| z \|_{E_{\varphi }}}} ) \leq \| y \|_{E} = \| \varphi(\frac{| z|}{\| z\|_{E_{\varphi}}})\|_E \leq 1.
\end{equation*}
Thus $\| z_1 \|_{E_{\varphi _{1}}}\leq \sqrt{D \| z \|_{E_{\varphi }}}$ and similarly 
$\| z_2 \|_{E_{\varphi _2}} \leq \sqrt{D \| z \|_{E_{\varphi }}}$. 
Since $| z | = z_1 z_2$ it follows that $z\in E_{\varphi_1}\odot E_{\varphi _2}$ and 
$\| z \|_{E_{\varphi_1}\odot E_{\varphi _2}}\leq D \| z \|_{E_{\varphi }}$.

$(a2)$ If $b_{\varphi }<\infty $ and $L^{\infty}\hookrightarrow E$, then $E_{\varphi } = L^{\infty}$ 
with equivalent norms and clearly $E_{\varphi } = L^{\infty }\hookrightarrow E_{\varphi
_1}\odot E_{\varphi _2}$. Suppose $b_{\varphi }=\infty$.
Set $v_0 =\varphi ^{-1}( u_0) $, where $u_0$ is from (\ref {right-ineq}) and let
$v >0$ be such that $\max [\varphi _1(v), \varphi _2(v)] \, \| \chi _{\Omega} \|_{E} \leq 1/2$.
For $\| z \|_{E_{\varphi }} = 1$ let $y = \varphi(|z|)$ and
\begin{equation*}
A = \{ t\in {\it supp} z: | z( t)| \geq v_0\}, \, B = {\it supp} z \backslash A = \{
t\in {\it supp} z: | z(t)| < v_0\},
\end{equation*}
Define 
\begin{equation*}
z_i ( t) = \left\{ 
\begin{array}{ccc}
 \sqrt{\frac{| z(t)| }{\varphi _{1}^{-1}( y(t)) \varphi _{2}^{-1} ( y( t) }}\,\varphi _{i}^{-1}( y( t)), & \text{if} & 
t\in A, \\ 
\sqrt{| z(t) | }, &  & t\in B, \\ 
0, &  & \text{otherwise},
\end{array}
\right.
\end{equation*}
for $i=1, 2$. Since $\varphi (v_0) >0$ the functions $z_i$ are well defined. 
If $t\in A$, then
\begin{equation*}
z_i( t) = \sqrt{\frac{| z(t)| }{\varphi _{1}^{-1}(y(t)) \varphi _{2}^{-1} (y(t))}} \, \varphi _{i}^{-1}( y(t))
\leq 
\sqrt{\frac{D | z(t) | }{\varphi^{-1}( y(t)) }} \,\varphi _{i}^{-1} (y( t)) \leq \sqrt{D}\, \varphi _{i}^{-1}( y(t)),
\end{equation*}
whence
\begin{equation*}
I_{\varphi _1} ( \frac{z_1}{2\sqrt{D}}\chi _{A}) \leq \frac{1}{2}I_{\varphi _1} 
( \frac{z_1}{\sqrt{D}}\chi _{A}) \leq \frac{1}{2} \| y \|_{E} \leq \frac{1}{2},
\end{equation*}
and 
\begin{equation*}
I_{\varphi _1} ( \frac{v z_{1}}{\sqrt{v_0}} \chi _{B}) = \| \varphi_1( \frac{v z_{1}}{\sqrt{v_0}} \chi _{B}) \|_E \leq
\varphi _1 (v) \| \chi _{\Omega} \|_{E} \leq \frac{1}{2}.
\end{equation*}
Then, for $\lambda = \max \{\frac{\sqrt{v_0}}{v}, 2 \sqrt{D} \}$, we obtain
\begin{equation*}
I_{\varphi _1}( \frac{z_1} {\lambda}) \leq I_{\varphi _1}( \frac{z_1} {\lambda}\, \chi _{A} ) + 
I_{\varphi _1}( \frac{z_1} {\lambda}\, \chi _{B}) \leq I_{\varphi _1} ( \frac{z_1}{2\sqrt{D}}\chi _{A}) + 
I_{\varphi _1} ( \frac{v z_{1}}{\sqrt{v_0}} \chi _{B}) \leq 1.
\end{equation*}
Thus $\| z_1 \|_{E_{\varphi _1}} \leq \lambda$ and similarly $\| z_2 \|_{E_{\varphi _2}} \leq \lambda$. 
Since $| z | = z_1\, z_2$ it follows that  $z\in E_{\varphi _1}\odot E_{\varphi _2}$ and 
$\| z \|_{E_{\varphi _1}\odot E_{\varphi _2}}\leq \lambda^2$.
Consequently $\| z \|_{E_{\varphi _1}\odot E_{\varphi _2}}\leq \lambda^2 \, \| z \|_{E_{\varphi }}$
for each $z\in E_{\varphi }$.
\vspace{2mm}

$(a3)$ Since $\varphi^{-1} \prec \varphi _{1}^{-1}\varphi _{2}^{-1}$ for small arguments it follows 
that for any $u_1 > u_0$ there is a constant $D_1 \geq D$ such that 
\begin{equation} \label{14}
\varphi ^{-1}(u) \leq D_{1} \varphi _{1}^{-1} (u)
\varphi_{2}^{-1} (u)  
\end{equation}
for any $u\leq u_{1}$. We follow the same way as in the proof of $(a1) $ replacing $D$ by $D_{1}$ 
from (\ref{14}) for $u_{1}= M ,$ where $M$ is the constant of the inclusion 
$E\overset{M}{\hookrightarrow }L^{\infty }$. 
\vspace{2mm}

$(b)$ Suppose that condition $\varphi^{-1} \prec \varphi _{1}^{-1}\varphi _{2}^{-1}$ for large arguments 
is not satisfied. Then there is a sequence $( u_n) $ with 
$u_n \nearrow \infty $ such that $2^{n}\varphi _{1}^{-1}( u_n) \varphi _{2}^{-1}(
u_n) \leq \varphi ^{-1}( u_n)$ for all  $n\in \mathbb N$.

We repeat a construction of  the sequence $( z_n) $, as it was given in [KLM12] in the proof of Theorem 4.2(i), 
showing that $\frac{\| z_{n}\|_{E_{\varphi _1} \odot E_{\varphi _2}}}{\| z_{n} \|_{E_{\varphi }}}\rightarrow \infty$ 
as $n \rightarrow \infty$.

Since $E_{a} \not= \{ 0\}$ it follows that there is a nonzero $0 \leq x\in E_{a}$ and so there is a set $A$ of positive 
measure such that $\chi _{A}\in E_{a}$. Of course, for large enough $n$ one has $\| u_{n}\chi_A \|_{E}\geq 1$. 
Applying Dobrakov's result from [Do74] we conclude that the submeasure $\omega ( B) = \| u_{n}\chi _B \|_E$ 
has the Darboux property. Consequently, for each $n\in \mathbb N$ there exists a set $A_n$ such that 
$\| u_n \chi_{A_n} \|_E = 1$. Define 
\begin{equation*}
x_n = \varphi _1^{-1}( u_n) \chi _{A_n}, ~ y_n = \varphi _2^{-1}( u_n) \chi _{A_n}  \text{ and } z_n = x_n\, y_n.
\end{equation*}
Let us consider two cases:

$1^0$. Let either $b_{\varphi_1} = \infty$ or $b_{\varphi_1} < \infty$ and $\varphi_1(b_{\varphi_1}) = \infty$. 
Then $\varphi_1(\varphi_1^{-1}(u)) = u$ for $u \geq 0$ and for $0 < \lambda < 1$, by the convexity of 
$\varphi_1$, we obtain
$$
I_{\varphi _1} ( \frac{x_{n}}{\lambda }) = \| \varphi_1 ( \frac{\varphi _1^{-1} ( u_n) }{\lambda })
\chi _{A_{n}} \|_E \geq \frac{1}{\lambda } \| \varphi_1 ( \varphi _1^{-1}( u_n)) \chi_{A_n} \|_E =
\frac{1}{\lambda } u_n \| \chi _{A_n} \|_E > 1.
$$

$2^0$. Let $b_{\varphi_1} < \infty$ and $\varphi_1(b_{\varphi_1}) < \infty$. Then, for sufficiently large $n$ 
and $0 < \lambda < 1$, we have $I_{\varphi_1} ( \frac{x_n}{\lambda }) =\infty $.

In both cases $\| x_n \|_{E_{\varphi_1}} \geq 1$ and similarly $\| y_n \|_{E_{\varphi_2}} \geq 1$. 
Applying Theorem 2 we get
\begin{eqnarray*}
\| z_n \|_{E_{\varphi _1}\odot E_{\varphi _2}}
&=&
\varphi _1^{-1}( u_n ) \varphi _2^{-1}( u_n )\, f_{E_{\varphi _1}\odot E_{\varphi _2}}( m(A _n))\\
&=&
\varphi _1^{-1}( u_n) \varphi _2^{-1}( u_n) \, f_{E_{\varphi _1}} ( m(A _n)) f_{E_{\varphi _2}}(m(A _n)) \\
&=&
\| x_n \|_{E_{\varphi _1}} \| y_n \|_{E_{\varphi _2}} \geq 1.
\end{eqnarray*}
On the other hand, using the relation between functions on a sequence $(u_n)$ and the fact that 
$\varphi(\varphi^{-1}(u)) \leq u$ for $u > 0$ we obtain
$$
I_{\varphi }( 2^n z_n ) = \| \varphi (2^n \varphi _1^{-1} ( u_n) \varphi _2^{-1} (u_n) ) \chi _{A_n} \|_E 
\leq \| \varphi ( \varphi ^{-1}( u_n) ) \chi _{A_n} \|_E \leq \| u_n \, \chi _{A_n} \|_E = 1, 
$$
i.e., $\| z_n \|_{E_{\varphi }} \leq 1/2^{n}$ which gives
$\frac{\| z_n \|_{E_{\varphi _1} \odot E_{\varphi _2}}}{\| z_n \|_{E_{\varphi }}} \geq 2^n \rightarrow \infty$ as 
$n \rightarrow \infty$.

$(c)$ It can be done by combining methods from the proof of
Theorem 5(b) and Theorem A(c).

$(d)$ Suppose that condition $\varphi^{-1} \prec \varphi _{1}^{-1}\varphi _{2}^{-1}$ for small arguments 
is not satisfied. Then we can find a sequence $u_{n}\rightarrow 0$ such that 
$2^n \varphi _1^{-1}( u_n) \varphi _2^{-1} (u_n) \leq \varphi^{-1} ( u_n) $ for all $n \in \mathbb N$.

Assumptions on a sequence space $e$ gives $\lim_{m \rightarrow \infty} \|\sum_{k=0}^m e_{k} \|_e = \infty$. 
Thus, for each $n \in \mathbb N$ there is a number $m_n$ such that 
\begin{equation*}
u_n \, \| \sum_{k=0}^{m_n} e_{k} \|_e \leq 1< u_n\, \| \sum_{k=0}^{m_n+1} e_k \|_e.
\end{equation*}
By symmetry of $e$, $\sup_{k \in \mathbb N} \| e_k \|_e = \| e_1\|_e = M$. Therefore 
$u_n \, \| \sum_{k=1}^{m_n} e_{k} \|_e \rightarrow 1$ as $n\rightarrow \infty $. Put 
$$
x_n = \varphi _1^{-1}( u_n) \, \sum_{k=1}^{m_n} e_k, ~ y_n = \varphi _2^{-1}( u_n) \, \sum_{k=1}^{m_n} e_k ~ 
{\rm and} ~ z_n = x_n \, y_n.
$$
Then $I_{\varphi _1} ( x_n) \leq u_n\, \| \sum_{k=1}^{m_n} e_{k} \|_e  \leq 1$ and 
\begin{equation*}
I_{\varphi _1}( 2x_{n}) = \varphi _1( 2\varphi_1^{-1}( u_n))\, \| \sum_{k=1}^{m_n} e_k \|_e \geq 2 u_n \, \|
\sum_{k=1}^{m_n} e_k \|_e \rightarrow 2 \text{ as } n \rightarrow \infty .
\end{equation*}
Therefore, for $n$ large enough $1\geq \| x_n \|_{e_{\varphi _1}}\geq 1/2$ as well as 
$1\geq \| y_n \|_{e_{\varphi _2}} \geq 1/2$. Consequently, explaining like in (b) one has 
$\| z_n \|_{E_{\varphi _1}\odot E_{\varphi _2}} \geq 1/4$ and $I_{\varphi } ( 2^n z_n) \leq 1$,
which gives $\frac{\| z_n \|_{E_{\varphi _1} \odot E_{\varphi _2}}}{\| z_n \|_{E_{\varphi }}} \geq 2^{n-2} \rightarrow \infty$ 
as $n \rightarrow \infty$ and the proof of Theorem 5 is complete.
\end{proof} 

To formulate results on equality of product spaces we need to introduce equivalences 
of inverses of Young functions $\varphi _1, \varphi _2$ and $ \varphi$.
The symbol $\varphi _{1}^{-1}\varphi _{2}^{-1}\approx \varphi ^{-1}$ for all
arguments [for large arguments] (for small arguments) means that $\varphi
_{1}^{-1}\varphi _{2}^{-1}\prec \varphi ^{-1}$ and $\varphi ^{-1}\prec
\varphi _{1}^{-1}\varphi _{2}^{-1},$ that is provided there are constants $C, D > 0$ 
[there are constants $C, D, u_0 >0$] (there are constants $C, D, u_0 >0$) such 
that the inequalities
\begin{equation} \label{equivalence}
C\varphi _{1}^{-1} ( u) \varphi _{2}^{-1} ( u) \leq \varphi ^{-1} ( u) \leq
D\varphi _{1}^{-1} ( u) \varphi_{2}^{-1} ( u)
\end{equation}
hold for all $u>0$ [for all $u\geq u_0$] (for all $0 < u\leq u_0$),
respectively.
\vspace{3mm}

From the above Theorem A and Theorem 5 we obtain immediately results on the product of 
Calder\'on-Lozanovski{\u \i} $E_{\varphi}$-spaces which are generalizations of the results known for 
Orlicz spaces.

\vspace{3mm} 
{\bf Corollary 6.} \label{Corollary6} {\it Let $\varphi_1, \varphi_2$ and $\varphi$ be three Young functions.
\vspace{-2mm}
\begin{itemize}
\item[$(a)$] Suppose that $E$ is a symmetric Banach function space with the Fatou property, 
$L^{\infty }\not\hookrightarrow E$ and ${\it supp} E_{a} = \Omega$. Then 
$E_{\varphi _1} \odot E_{\varphi _2} = E_{\varphi } $ if and only if
$\varphi _{1}^{-1}\varphi _{2}^{-1}\approx \varphi ^{-1}$ for all arguments.

\item[$(b)$] Suppose that $E$ is a symmetric Banach function space with the Fatou property, $L^{\infty }\hookrightarrow E$ 
and $E_{a} \neq \{ 0\} $. Then $E_{\varphi _1} \odot E_{\varphi _2} = E_{\varphi } $ if and 
only if $\varphi _{1}^{-1}\varphi _{2}^{-1}\approx \varphi ^{-1}$ for large arguments.

\item[$(c)$] Suppose that $e$ is a symmetric Banach sequence space with the Fatou property and 
order continuous norm. Then $e_{\varphi _1} \odot e_{\varphi _2} = e_{\varphi } $
if and only if $\varphi _{1}^{-1}\varphi _{2}^{-1}\approx \varphi ^{-1}$ for small arguments.
\end{itemize}
}

The following construction appeared in [ZR67] and in [DK67]: for two Young functions $\varphi _1,\varphi _2$ 
(or even for only the so-called $\varphi$-functions) one can define a new function 
$\varphi _1 \oplus \varphi _2$ by the formula
\begin{equation}
( \varphi _1 \oplus \varphi _2) (u) = \inf_{u = vw } [ \varphi _1(v) +\varphi _2( w)] 
= \inf_{v > 0} \, [\varphi_1(v) + \varphi_2(\frac{u}{v})],
\end{equation}
for $u\geq 0$. This operation was investigated in [ZR67], [BO78], [Ma89], [MP89] and [St92]. Note that 
$\varphi _1 \oplus \varphi _2$ is non-decreasing, left-continous function and is 0 at $u = 0$. 
Moreover, the proof of the following estimates can be found in [ZR67, pp. 267, 271] and [St92, Theorem 1]: 
$$
\varphi^{-1}( t) \leq \varphi _1^{-1}( t) \varphi_2^{-1}( t) \leq \varphi ^{-1} ( 2 t) ~~ {\rm for ~any} ~ t > 0,
$$
where $\varphi _1, \varphi _2$ are nondegenerate Orlicz functions (the proof in general case is not difficult) 
$\varphi$-functions and $\varphi =  \varphi _1 \oplus \varphi _2$. 
The function $\varphi$ need not be convex even if both $\varphi _1$ and $\varphi _2$ are 
convex functions. However, if $\varphi$ is a convex function, then 
$\varphi^{-1} \leq \varphi_1^{-1} \varphi_2^{-1} \leq 2 \varphi^{-1}$ and, by Theorem A(a1) and 5(a1), we obtain 
$E_{\varphi_1 } \odot E_{\varphi_2} = E_{\varphi }$.

We will prove the last result without explicit assumption that $\varphi$ is convex, but to do this we need 
to extend definition of the Calder\'{o}n-Lozanovski{\u \i} $E_{\varphi}$-space in this case (cf.  [KMP03] 
for definition and some results).

For a non-decreasing and left-continous function $\varphi: [0, \infty) \rightarrow [0, \infty)$ with 
$\varphi(0) = 0$ assume that there exist $C, \alpha > 0$ such that
\begin{equation} \label{equation18}
\varphi(st) \leq C\, t^{\alpha}\, \varphi(s) ~~ {\rm for ~all} ~ s > 0 ~ {\rm and} ~ 0 < t < 1.
\end{equation}
Then the Calder\'{o}n-Lozanovski{\u \i} $E_{\varphi}$-space is a quasi-Banach ideal space (since $\varphi $ need 
not be convex) with the quasi-norm
$$
\| x\|_{E_{\varphi}} = \inf \{\lambda > 0: \| \varphi(|x|/\lambda) \|_E \leq 1\}.
$$
Note that for a convex function $\varphi $ the condition (\ref{equation18}) holds with $C = \alpha = 1$ and the space 
$E_{\varphi }$ is normable when the condition (\ref{equation18}) holds with $\alpha \geq 1$.

We know that if $\varphi_1$ and $\varphi_2$ are convex functions, then $\varphi = \varphi _1 \oplus \varphi _2$ is 
not necessary a convex function but the condition (\ref{equation18}) holds with $C = 1$ and $\alpha= 1/2$ and the 
space  $E_{\varphi }$ is a quasi-Banach ideal space. In fact, for any $s > 0$ and $0 < t < 1$ we have
\begin{eqnarray*}
\varphi(st) 
&=& 
\inf_{st = v w} [\varphi _1(v) +\varphi _2( w)] =
\inf_{s = ab} [\varphi _1(\sqrt{t} \, a) +\varphi _2( \sqrt{t}\, b)] \\
&\leq&
\sqrt{t} \, \inf_{s = ab} [\varphi _1(a) +\varphi _2(b)] = \sqrt{t} \, \varphi(s).
\end{eqnarray*}

{\bf Theorem 6.} {\it Let $\varphi _{1},\varphi _{2}$ be Young functions. If 
$\varphi: = \varphi _1 \oplus \varphi _2$, then $E_{\varphi _1}\odot
E_{\varphi _2} = E_{\varphi }$, where $E$ is a Banach ideal space with the Fatou property.
}

\begin{proof}  We prove that $E_{\varphi _1} \odot E_{\varphi _2}\hookrightarrow
E_{\varphi }$. By definition of $\varphi _1 \oplus \varphi _2$ one has
\begin{equation*}
\varphi ( uv) \leq \varphi _1( u) +\varphi_2(v)
\end{equation*}
for each $u, v > 0$ with $\varphi _1( u) < \infty, \varphi _2(v) < \infty $. Let 
$z \in E_{\varphi _1}\odot E_{\varphi _2}, z\neq 0$, and take arbitrary $0 \leq x \in E_{\varphi _1}, 0 \leq y \in E_{\varphi _2}$ 
with $| z | = x y$. Since the condition (\ref{equation18}) holds with $C = 1$ and $\alpha= 1/2$ it follows that for 
$0 < t < 1/4$ we obtain
\begin{eqnarray*}
I_{\varphi }( \frac{t\, z}{\| x \|_{E_{\varphi_1}} \| y \|_{E_{\varphi _2}}}) 
&\leq&
\sqrt{t} \, \| \varphi ( \frac{x}{\| x \|_{E_{\varphi _1}}} \frac{y}{\| y \|_{E_{\varphi _2}}} ) \|_E \\
&\leq&
\sqrt{t} \, \left[ \| \varphi _1 ( \frac{x}{\| x \|_{E_{\varphi _1}}} ) \|_E + 
\| \varphi_2 ( \frac{y}{\| y \|_{E_{\varphi _2}}} ) \|_E \right ] \\
&\leq&
2\, \sqrt{t} < 1.
\end{eqnarray*}
Thus $\| z \| _{E_{\varphi }} \leq \frac{1}{t} \, \| x \|_{E_{\varphi _1}} \| y \|_{E_{\varphi _2}}$
and consequently $\| z \|_{E_{\varphi }} \leq \frac{1}{t} \| z \|_{E_{\varphi _1} \odot E_{\varphi _2}}$. 
Thus $E_{\varphi _1} \odot E_{\varphi _2} \overset{1/t}{\hookrightarrow } E_{\varphi }$.

The proof of the imebedding $E_{\varphi }\hookrightarrow E_{\varphi _1}\odot
E_{\varphi _2}$ is exactly the same as the proof of Theorem 5(a1) since the convexity of 
$\varphi $ has not been used there, which proves the theorem.
\end{proof} 

\vspace{3mm}
\begin{center}
\textbf{5. The product of Lorentz and Marcinkiewicz spaces}
\end{center}

Before proving results on the product of Lorentz and Marcinkiewicz spaces on $I = (0, 1)$ 
or $I = (0, \infty)$ we need some auxiliary lemmas on the Calder\'on construction and notion 
of the dilation operator.

The {\it dilation operator} $D_s, s > 0,$ defined by $D_sx(t) = x(t/s) \chi_I(t/s), t \in I$ is bounded 
in any symmetric space $E$ on ${\rm I}$ and $\| D_s\|_{E \rightarrow E} \leq \max(1, s)$ (see 
[Sh68, Lemma 1] in the case $I = (0, 1)$, [KPS82, pp. 96-98] for $I = (0, \infty)$ and [LT79, p. 130] 
for both cases). Moreover, the {\it Boyd indices} of $E$ are defined by
\begin{equation*}
\alpha_E = \lim_{s \rightarrow 0^+} \frac{\ln \| D_s\|_{E \rightarrow E} }{\ln s}, \beta_E = 
\lim_{s \rightarrow \infty} \frac{\ln \| D_s\|_{E \rightarrow E} }{\ln s},  
\end{equation*}
and we have $0 \leq \alpha_E \leq \beta_E \leq 1$.

\vspace{3mm}
{\bf Lemma 3.} \label{Lemma3} {\it Let $E, F$ be symmetric function spaces on $I$ and $0 < \theta < 1$. Then
\begin{equation*}
\frac{1}{2} \, \| z^* \|_{E^\theta F^{1-\theta}}^* \leq  \| z^* \|_{E^\theta F^{1-\theta}} \leq  \| z^* \|_{E^\theta F^{1-\theta}}^*,
\end{equation*}
where  $\| z^* \|_{E^\theta F^{1-\theta}}^* : = \inf \{ \max ( \| x^* \|_E, \| y^* \|_F): z^* \leq (x^*)^\theta \, (y^*)^{1-\theta} , 
x \in E_+, y \in F_+ \}.$}

\begin{proof} Since
\begin{equation*}
\| x^* \|_E = \| D_2 D_{\frac{1}{2}} x^* \|_E \leq \| D_2\|_{E \rightarrow E} \, \| D_{\frac{1}{2}} x^* \|_E \leq 2\, 
\| D_{\frac{1}{2}} x^* \|_E.
\end{equation*}
and
$(| x|^{\theta} \, |y|^{1-\theta})^*(t) \leq x^*(t/2)^{\theta} \, y^*(t/2)^{1-\theta}$ for any $t \in I$ (cf. [KPS82, p. 67]) 
it follows that
\begin{eqnarray*}
\| z^* \|_{E^\theta F^{1-\theta}}
&=& 
\inf \{ \max ( \| x \|_E, \| y \|_F): z^* \leq x^\theta \, y^{1-\theta} , x \in E_+, y \in F_+ \} \\
&\geq&
\inf \{ \max ( \| x \|_E, \| y \|_F): z^*(t) \leq x^*(t/2)^\theta \, y^*(t/2)^{1-\theta} , x \in E_+, y \in F_+ \}\\
&=&
\inf \{ \max ( \| D_{\frac{1}{2}} x^* \|_E, \| D_{\frac{1}{2}} y^* \|_F): z^*(t) \leq x^*(t)^\theta \, y^*(t)^{1-\theta} , x \in E_+, 
y \in F_+ \}\\
&\geq&
\frac{1}{2} \, \inf \{ \max ( \| x^* \|_E, \| y^* \|_F): z^*(t) \leq x^*(t)^\theta \, y^*(t)^{1-\theta} , x \in E_+, y \in F_+ \} \\
&=&
\frac{1}{2} \, \| z^* \|_{E^\theta F^{1-\theta}}^*.
\end{eqnarray*}
The other estimate is clear and the lemma follows.
\end{proof} 

As a consequence of  representation (\ref{equation8}) and the above lemma with $\theta = 1/2$ we obtain

\vspace{3mm} 
{\bf Corollary 7.} \label{Corollary3} {\it Let $E, F$ be symmetric function spaces on $I$. Then
\begin{equation*}
\| z^*\|_{E\odot F}\leq \inf \{ \| x \|_E \, \| y \|_F: z^* \leq x^* \, y^*, x \in E_+, y \in F_+ \} \leq 2\, \| z^*\|_{E\odot F}.
\end{equation*}

}

The idea of the proof of the next result is coming from Calder\'on [Ca64, Part 13.5]. For a Banach function 
space $E$ on $I = (0, 1)$ or $(0, \infty)$ define new spaces (symmetrizations of $E$) $E^{(*)}$ and $E^{(**)}$ as
$$
E^{(*)} = \{x \in L^0(I): x^* \in E \}, ~~ E^{(**)}= \{x \in L^0(I): x^{**} \in E \}
$$
with the functionals $\| x\|_{E^{(*)}} = \| x^* \|_E$ and $\| x\|_{E^{(**)}} = \| x^{**}\|_E$. If $C_E$ denotes 
the smallest constant $1 \leq C < \infty$ such that 
\begin{equation} \label{D2}
\| D_2 x^* \|_E \leq C \, \| x^*\|_E ~~ {\rm for ~ all} ~  x^* \in E,
\end{equation}
then $E^{(*)}$ is a quasi-Banach symmetric space. The space $E^{(**)}$ is always a Banach symmetric space. 
Consider the Hardy operator $H$ and its dual $H^*$ defined by
\begin{equation}\label{Hardy}
Hx(t) = \frac{1}{t} \int_0^t x(s)\, ds, ~ H^*x(t) = \int_t^l \frac{x(s)}{s}\, ds ~ {\rm with} ~ l = m(I), ~ t \in I.
\end{equation}

{\bf Remark 7.} If $E$ is a Banach function space on $I$ and operator $H$ is bounded in $E$, then 
(\ref{D2}) holds with $C_E \leq 2\, \| H\|_{E \rightarrow E}$. This follows directly from the estimates
$$
\| Hx^*\|_E = \| \int_0^1 x^*(st)\, ds \|_E \geq  \| \int_{0}^{1/2} x^*(st)\, ds \|_E \geq \frac{1}{2} \, \| x^*(t/2)\|_E.
$$

As we already mentioned before the Calder\'on spaces $E^{\theta} F^{1-\theta}$ can 
be also defined for quasi-Banach spaces $E, F$ (cf. [Ov82], [Ni85], [KMP03]).

\vspace{3mm}
{\bf Lemma 4.} \label{Lemma4} {\it Let $E$ and $F$ be Banach function spaces on $I$ and 
$0 < \theta < 1$. Suppose that both operators $H, H^*$ are bounded in $E$ and $F$. Then
\begin{equation} \label{embedding21}
(E^{(*)})^{\theta}(F^{(*)})^{1-\theta} \overset{C_1}{\hookrightarrow } (E^{\theta} F^{1-\theta})^{(*)} \overset{C_2}{\hookrightarrow } (E^{(*)})^{\theta}(F^{(*)})^{1-\theta},
\end{equation}
where $C_1 = C_E^{\theta} C_F^{1-\theta}, C_2 = \| H H^* \|_{E \rightarrow E}^{\theta} \| H  H^* \|_{F \rightarrow F}^{1- \theta}$ and
\begin{equation}  \label{embedding22}
(E^{(**)})^{\theta}(F^{(**)})^{1-\theta} \overset{1}{\hookrightarrow } (E^{\theta} F^{1-\theta})^{(**)} \overset{C_3}{\hookrightarrow } (E^{(**)})^{\theta}(F^{(**)})^{1-\theta},
\end{equation}
where $C_3 = [ \| H \|_{E \rightarrow E}\,  \| H H H^* \|_{E \rightarrow E}]^{\theta}\, [ \| H \|_{F \rightarrow F} \, \| H H H^* \|_{F \rightarrow F}]^{1-\theta}.$
}

\begin{proof} Embeddings (\ref{embedding21}). Let $z \in (E^{(*)})^{\theta}(F^{(*)})^{1-\theta}$. Then 
$|z| \leq \lambda |x|^{\theta} |y|^{1-\theta}$ for 
some $\lambda > 0$ and $\| x^*\|_E \leq 1, \| y^* \|_F \leq 1$. Thus 
$$
z^*(t) \leq \lambda(|x|^{\theta} |y|^{1-\theta})^*(t) \leq \lambda x^*(t/2)^{\theta} y^*(t/2)^{1-\theta} 
= \lambda C_1 (\frac{x^*(t/2)}{C_E})^{\theta} (\frac{y^*(t/2)}{C_F})^{1-\theta} 
$$
for any $t \in I$, which means that $z \in (E^{\theta} F^{1-\theta})^{(*)}$ with the norm $\leq \lambda C_1$.
\vspace{2mm}

On the other hand, if $z \in (E^{\theta} F^{1-\theta})^{(*)}$, then $z^* \in E^{\theta} F^{1-\theta}$ and so
$$
z^* \leq \lambda |x|^{\theta} |y|^{1-\theta} ~ {\rm with ~some} ~\lambda > 0, \|x\|_E \leq 1, \| y\|_F \leq 1.
$$
The following equality is true
\begin{equation} \label{equality23}
HH^*x(t) = Hx(t) + H^*x(t), ~ t\in I.
\end{equation}
In fact, using the Fubini theorem, we obtain for $x \geq 0$
\begin{eqnarray*}
HH^*x(t) 
&=& 
\frac{1}{t} \int_0^t (\int_s^l \frac{x(r)}{r}\, dr)\, ds = \frac{1}{t} \int_0^t (\int_0^r ds)\, \frac{x(r)}{r}\, dr\\
&+& 
\frac{1}{t} \int_t^l (\int_0^t ds) \, \frac{x(r)}{r}\, dr = Hx(t) + H^*x(t).
\end{eqnarray*}
Using then equality (\ref{equality23}) and twice H\"older-Rogers inequality we obtain
\begin{eqnarray*}
z^* 
&\leq& 
H(z^*) \leq H(z^*) + H^*(z^*) = H H^*(z^*)\\
&\leq&
\lambda HH^*(|x|^{\theta} |y|^{1-\theta}) \leq \lambda H[ (H^*|x|)^{\theta} (H^*|y|)^{1-\theta}]\\
&\leq&
\lambda [HH^*(|x|)]^{\theta} \, [HH^*(|y|)]^{1-\theta}.
\end{eqnarray*}
By the Ryff theorem there exists a measure-preserving transformation $\omega: I \rightarrow I$ 
such that $|z| = z^*(\omega)$ a.e. (cf. [BS88], Theorem 7.5 for $I = (0, 1)$ or Corollary 7.6 for $I = (0, \infty)$ 
under additional assumption that $z^*(\infty) =0$). Thus
\begin{equation*}
|z| = z^*(\omega) \leq \lambda [HH^*(|x|)(\omega)]^{\theta} \, [HH^*(|y|)(\omega)]^{1-\theta} = \lambda u^{\theta} v^{1-\theta}.
\end{equation*}
Since $H^*|x|$ is non-increasing function it follows that $HH^*|x|$ is also non-increasing function and
$HH^*|x| = [HH^*|x|]^* = [(HH^*|x|)(\omega)]^*$. Similarly with $HH^*|y|$. Hence 
\begin{eqnarray*}
\| u \|_{E^{(*)}} 
&=& 
\| u^*\|_E = \| [(HH^*|x|)(\omega)]^*\|_E \\
&=&
 \| HH^*|x| \|_E \leq \| H H^*\|_{E \rightarrow E} \, \| x\|_E \leq \| H H^*\|_{E \rightarrow E} 
\end{eqnarray*}
and
\begin{eqnarray*}
\| v \|_{F^{(*)}} 
&=& 
\| v^*\|_F = \| [(HH^*|y|)(\omega)]^*\|_F\\
&=&
 \| HH^*|y| \|_F \leq \| H H^*\|_{F \rightarrow F} \, \| y\|_F \leq \| H H^*\|_{F \rightarrow F} ,
\end{eqnarray*}
which means that $z \in (E^{(*)})^{\theta}(F^{(*)})^{1-\theta}$ with the norm $\leq \lambda C_2$.

To finish the proof in the case $I = (0, \infty)$ we need to show $z^*(\infty) = 0$. If we will have 
 $z^*(\infty) = a > 0$, then $\lambda |x(t)|^{\theta} |y(t)|^{1-\theta} \geq a$ for almost all $t > 0$ and 
considering the sets $A = \{t > 0: |x(t)| \geq a/\lambda \}, ~B = \{t > 0: |y(t)| \geq a/\lambda\}$ we obtain 
$A \cup B = (0, \infty)$ up to the set of measure zero. Then
$$
H^*|x|(t) = \int_t^{\infty} \frac{|x(s)|}{s}\, ds \geq \int_{A \cap(t, \infty)} \frac{a}{\lambda s} \, ds
$$
and
$$
H^*|y|(t) = \int_t^{\infty} \frac{|y(s)|}{s}\, ds \geq \int_{B \cap(t, \infty)} \frac{a}{\lambda s} \, ds,
$$
which means $H^*|x|(t) + H^*|y|(t) = + \infty$ for all $t > 0$. Since
$$
(0, \infty) = \{ t > 0: H^*|x|(t) = \infty \} \cup \{ t > 0: H^*|y|(t) = \infty \}
$$ 
(maybe except the set of measure zero) it follows that $H^*|x| \notin E$ or $H^*|y| \notin F$, which is 
a contradiction.
\vspace{2mm}

Embeddings (\ref{embedding22}). Let $z \in (E^{(**)})^{\theta}(F^{(**)})^{1-\theta}$. Then 
$|z| \leq \lambda |x|^{\theta} |y|^{1-\theta}$ for some $\lambda > 0$ and $\| x^{**}\|_E \leq 1, \| y^{**} \|_F \leq 1$. 
Thus 
$$
z^{**}(t) \leq \lambda(|x|^{\theta} |y|^{1-\theta})^{**}(t) \leq \lambda x^{**}(t)^{\theta} y^{**}(t)^{1-\theta}
$$
for any $t \in I$, which means that $z \in (E^{\theta} F^{1-\theta})^{(**)}$ with the norm $\leq \lambda$.

On the other hand, if $z \in (E^{\theta} F^{1-\theta})^{(**)}$, then $z^{**} \in E^{\theta} F^{1-\theta}$ and 
repeating the above arguments we obtain
\begin{eqnarray*}
|z| 
&=& 
z^*(\omega) \leq z^{**}(\omega) = Hz^*(\omega) \\
&\leq&  
\lambda \, [ H H H^* |x|(\omega)]^{\theta} \, [H H H^*(|y|)(\omega)]^{1-\theta} = \lambda \, u_1^{\theta} v_1^{1-\theta}.
\end{eqnarray*}
Since $H H H^*|x| = [H H H^*|x|]^* = [(H H H^*|x|)(\omega)]^*$ it follows that
\begin{eqnarray*}
\| u_1\|_{E^{(**)}} 
&=& 
\| u_1^{**} \|_E = \| H u_1^* \|_E \leq  \| H\|_{E \rightarrow E} \, \| u_1^* \|_E\\
&=& 
\| H\|_{E \rightarrow E} \, \| [(H H H^*|x|)(\omega)]^* \|_E =  \| H\|_{E \rightarrow E} \, \| H H H^*|x| \|_E \\
&\leq&
\| H\|_{E \rightarrow E} \,  \| H H H^* \|_{E \rightarrow E}  \| x\|_E \leq \| H\|_{E \rightarrow E} \,  \| H H H^* \|_{E \rightarrow E}
\end{eqnarray*}
and
\begin{eqnarray*}
\| v_1\|_{F^{(**)}} 
&=& 
\| v_1^{**} \|_F = \| H v_1^* \|_F \leq  \| H\|_{F \rightarrow F} \, \| v_1^* \|_F\\
&=& 
\| H\|_{F \rightarrow F} \, \| [(H H H^*|y|)(\omega)]^* \|_F =  \| H\|_{F \rightarrow F} \, \| H H H^*|y| \|_F \\
&\leq&
\| H\|_{F \rightarrow F} \,  \| H H H^* \|_{F \rightarrow F}  \| y\|_F \leq \| H\|_{F \rightarrow F} \,  \| H H H^* \|_{F \rightarrow F},
\end{eqnarray*}
which implies that $z \in (E^{(**)})^{\theta}(F^{(**)})^{1-\theta}$ with the norm $\leq \lambda C_3$, and the lemma follows.
\end{proof} 

Note that our proofs are working for both cases $I = (0, 1)$ and $I = (0, \infty)$. Our inclusions (\ref{embedding21}) were 
proved by Calder\'on but his result is true only in the case when $I = (0, \infty)$ (cf. [Ca64], pp. 167-169). 
Since he was working with the other composition $H^*H$, which in the case $I = (0, \infty)$ gives the equality 
$H^* H = H + H^*$. For $I = (0, 1)$ one gets another formula $H^* H x (t) = Hx(t) + H^*x(t) - \int_0^1 x(s)\, ds$, 
which not allows then to proof the same result in this case.

For the identification of product spaces we will need result on the Calder\'on construction for weighted Lebesgue 
spaces $L^p(w) = \{ x \in L^0(\mu): x w \in L^p(\mu)\}$ with the norm $\| x\|_{L^p(w)} = \| x w\|_{L^p}$, where 
$ 1 \leq p \leq \infty$ and $w \geq 0$ . Result for $p_0 = p_1$ was given in [Ov84, p. 459] and for general Banach 
ideal spaces in [KM03, Theorem 2] (for $1 \leq p_0, p_1 < \infty$ it also follows implicitely from [BL76, Theorem 5.5.3] 
and results on relation between the complex method and the Calder\'on  construction). We present here a direct proof.

\vspace{3mm}
\textbf{Lemma 5.} \label{Lemma5} {\it Let $1 \leq p_0, p_1 \leq \infty$ and $0 < \theta < 1$. Then
\begin{equation} \label{equality24}
L^{p_0}(w_0)^{1-\theta} \, L^{p_1}(w_1)^{\theta} \equiv L^p(w),
\end{equation}
where $\frac{1}{p} = \frac{1-\theta}{p_0} + \frac{\theta}{p_1}$ and $w = w_0^{1-\theta} w_1^{\theta}$.
}

\begin{proof} Suppose $1 \leq p_0, p_1 < \infty$. If $x \in L^{p_0}(w_0)^{1-\theta} \, L^{p_1}(w_1)^{\theta}$, then
$|x| \leq \lambda |x_0|^{1-\theta} |x_1|^{\theta}$ with $\|x_0\|_{L^{p_0}(w_0)} \leq 1$ and $\|x_1\|_{L^{p_1}(w_1)} \leq 1$. Using the H\"older-Rogers inequality we obtain
\begin{eqnarray*}
\int |x\, w|^p \, d\mu 
&\leq& 
\lambda^p \int |x_0 w_0|^{(1-\theta)p} |x_1 w_1|^{\theta p} \, d\mu\\
&\leq& 
\lambda^p (\int |x_0 w_0|^{p_0} \, d\mu)^{(1-\theta)p/{p_0}}\, (\int |x_1 w_1|^{p_1}\, d\mu)^{\theta p/{p_1}} \\
&=&
\lambda^p \, \| x_0 w_0\|_{L^{p_0}}^{(1-\theta)p} \,  \| x_1 w_1\|_{L^{p_1}}^{\theta p},
\end{eqnarray*}
that is
$$
\| x\|_{L^p(w)} \leq \lambda \,  \| x_0 w_0\|_{L^{p_0}}^{1-\theta} \,  \| x_1 w_1\|_{L^{p_1}}^{\theta} \leq \lambda
$$
and so $x \in L^p(w)$ with the norm $\leq \lambda$.

\noindent
On the other hand, if $0 \neq x \in L^p(w)$ then, considering
$x_i(t) = \frac{|x(t) \,w(t)|^{p/{p_i}}}{\| x \|_{L^p(w)}^{p/{p_i}}} \frac{1}{w_i(t)}$ on the support of $w_i$ and $0$ otherwise ($i = 0, 1$), we obtain $|x(t)| = \| x\|_{L^p(w)} \, x_0(t)^{1-\theta}\, x_1(t)^{\theta}$ and
$$
\| x_i \|_{L^{p_i}(w_i)}^{p_i} = \int [x_i(t)\, w_i(t)]^{p_i}\, d\mu =\int \frac{|x(t) w(t)|^p}{\| x \|_{L^p(w)}^p} \, d\mu = 1.
$$
Therefore, $x \in L^{p_0}(w_0)^{1-\theta} \, L^{p_1}(w_1)^{\theta} $ with norm $\leq  \| x\|_{L^p(w)}$ and equality (\ref{equality24}) is proved. The proof for the case when one or both $p_0, p_1$ are $\infty$ is even simpler.
\end{proof} 

We want to calculate product spaces of Lorentz space $\Lambda_\phi$ and Marcinkiewicz space $M_\phi$ on $I$, 
where $\phi$ is a quasi-concave function on $I$ with $\phi(0^+) = 0$. We will do this, in fact, for some other closely 
connected spaces. Consider the Lorentz space $\Lambda_{\phi, 1}$ and more general Lorentz space $\Lambda_{\phi, p}$ 
with $0 < p < \infty$ on $I$ defined, respectively, as
\begin{equation*}
\Lambda_{\phi, 1} = \{ x\in L^{0}( I): \| x \|_{\Lambda_{\phi, 1}} = \int_I x^*(t) \frac{\phi(t)}{t}\, dt < \infty\},
\end{equation*}
\begin{equation*}
\Lambda_{\phi, p} = \{ x\in L^{0}( I): \| x \|_{\Lambda_{\phi, p}} = \left( \int_I [\phi(t)\, x^*(t)]^p \frac{dt}{t} \right)^{1/p}< \infty\},
\end{equation*}
Space $\Lambda_{\phi, 1}$ is a Banach space and if $\phi(t) \leq a t \phi^{\prime}(t)$ for all $t \in I$, then $\Lambda_{\phi, 1} 
\overset{1}{\hookrightarrow } \Lambda_{\phi} \overset{a}{\hookrightarrow } \Lambda_{\phi, 1}$ (space $\Lambda_{\phi} $ was defined in Part 1).
Consider also another Marcinkiewicz space $M_{\phi }^*$ than the space $M_{\phi}$ defined in Part 1, as 
\begin{equation*}
M_{\phi }^* = M_{\phi }^*(I) = \{ x\in L^{0}( I): \| x \|_{M_{\phi }^*} = \sup_{t\in I} \phi ( t) x^*(t) <\infty\}.
\end{equation*}
This Marcinkiewicz space need not be a Banach space and always we have $M_{\phi } \overset{1}{\hookrightarrow }M_{\phi }^*$. Moreover,  $M_{\phi }^* \overset{C}{\hookrightarrow } M_{\phi }$ if and only if
\begin{equation} \label{estimate25}
\int_0^t \frac{1}{\phi(s)} \, ds \leq C \frac{t}{\phi(t)} ~ {\rm for ~all} ~~ t \in I.
\end{equation}
In fact, since $\frac{1}{\phi} \in M_{\phi }^*$ then estimate (\ref{estimate25}) is necessary for the imbedding. On the other hand, if (\ref{estimate25}) holds and $x \in M_{\phi }^*$, then
\begin{eqnarray*}
\| x \|_{M_{\phi }} 
&=& 
\sup_{t \in I} \phi(t) x^{**}(t) = 
\sup_{t \in I} \frac{\phi(t)}{t} \int_0^t \frac{1}{\phi(s)}\, \phi(s) x^*(s) \, ds\\
&\leq&
\sup_{s \in I} \phi(s) x^*(s) \, \sup_{t \in I} \frac{\phi(t)}{t} \int_0^t \frac{1}{\phi(s)} \, ds \leq C \| x \|_{M_{\phi }^*}.
\end{eqnarray*}

We can consider spaces $\Lambda_{w, 1}, \Lambda_{w, p}$ and $M_w^*$ for more general weights 
$w \geq 0$, but then the problem of being quasi-Banach space or Banach space will appear. Such 
investigations can be found in [CKMP] and [KM04].

Since indices of the quasi-concave function on $I$ are useful in the formulation of further results let us 
define them. The lower index $p_{\phi, I}$ and upper index $q_{\phi, I }$ 
of a function $\phi $ on $I$ are numbers defined as 
\begin{equation*}
p_{\phi, I} = \lim_{t \rightarrow 0^{+}} \frac{\ln m_{\phi, I}( t) }{\ln t}, ~ q_{\phi, I} = \lim_{t \rightarrow \infty } \frac{\ln m_{\phi, I}(t) }{\ln t}, ~ {\rm where} ~~ m_{\phi, I}( t) = \sup_{s \in I, st \in I }\frac{\phi( st ) }{\phi( s) }. 
\end{equation*}
It is known (see, for example, [KPS82] and [Ma85], [Ma89]) that for a quasi-concave function $\phi $ on $[0, \infty)$ 
we have $0 \leq p_{\phi, [0, \infty)} \leq p_{\phi, [0, 1]} \leq q_{\phi, [0, 1]} \leq q_{\phi, [0, \infty)} \leq 1$. Moreover, 
estimate (\ref{estimate25}) is equivalent to $q_{\phi, I} < 1$.
We also need for a differentiable increasing function $\phi$ on $I$ with $\phi(0^+) = 0$ the {\it Simonenko indices} 
\begin{equation*}
s_{\phi, I} = \inf_{t \in I} \frac{t \phi^{\prime}(t)}{\phi(t)}, ~  \sigma_{\phi, I} = \sup_{t \in I} \frac{t \phi^{\prime}(t)}{\phi(t)}.
\end{equation*}
They satisfy $0 \leq s_{\phi, I}  \leq p_{\phi, I} \leq q_{\phi, I} \leq \sigma_{\phi, I} $ (cf.  [Ma85, p. 22] 
and [Ma89, Theorem 11.11]).

\vspace{3mm}
{\bf Theorem 7.} 
\vspace{-3mm}
\begin{itemize}
\item[$(i)$] {\it  If $\phi ,\psi $ are quasi-concave functions on $I$, then 
$ M_{\phi \psi }^*  \overset{1}{\hookrightarrow } M_{\phi }^* \odot M_{\psi }^*  \overset{2}{\hookrightarrow } 
M_{\phi \psi }^*$.
}

\item[$(ii)$]  {\it Let $\phi ,\psi $ and $\phi \psi$ be increasing concave functions on $I$ with 
$\phi( 0^+) = \psi( 0^+) =0$. If $s_{\phi, I} \geq a > 0$ and $s_{\phi \psi, I} \geq b > 0$, then
$\Lambda_{\phi}  \odot M_{\psi }^* \overset{4+4/a}\hookrightarrow  \Lambda_{\phi \psi } \overset{2 /b}\hookrightarrow \Lambda_{\phi}  \odot M_{\psi }^*$.}

\item[$(iii)$] {\it Let $\phi ,\psi $ be quasi-concave functions on $I$ such that $0 < p_{\phi, I}  \leq q_{\phi, I} < 1$ 
and $0 < p_{\psi, I}  \leq q_{\psi, I} < 1$, then 
\begin{equation*}
\Lambda_{\phi, 1 } \odot \Lambda _{\psi, 1 } = \Lambda _{\phi \psi, 1/2 }, ~ \Lambda_{\phi, 1 } \odot M _{\psi }^* = \Lambda _{\phi \psi, 1 }, ~M_{\phi}^* \odot M_{\psi}^* = M _{\phi \psi}^*
\end{equation*}
\vspace{-3mm}
with equivalent quasi-norms.}
\end{itemize}

\begin{proof} (i) For each $z \in M_{\phi \psi }^{\ast }$ one has $z^* \leq \frac{\| z \|_{M_{\phi \psi }^{\ast }}}{\phi \psi }$, 
but since $ \frac{\| z \|_{M_{\phi \psi }^{\ast }}}{\phi \psi } \in M_{\phi}^* \odot M_{\psi }^*$ it follows 
that $z \in M_{\phi }^* \odot M_{\psi }^*$ and $ \| z \|_{M_{\phi }^* \odot M_{\psi }^*} \leq 
\| z \|_{M_{\phi \psi}^*}$. 
\vspace{1mm}

If $z \in M_{\phi }^* \odot M_{\psi }^*$, then, by Corollary 7, we have $z^* \leq x^* y^*$ for some 
$x^* \in M_{\phi }^*, y^* \in M_{\psi }^*$ and 
$\inf \{ \| x\|_{M_{\phi }^*} \| y \|_{M_{\psi }^*}: z^* \leq x^* y^*\} \leq 2 \| z\|_{ M_{\phi }^* \odot M_{\psi }^*}$.
But $x^* y^* \leq \frac{\| x\|_{M_{\phi }^*}}{\phi }\, \frac{\| y\|_{M_{\psi }^*}}{\psi }$ and so
$$
\| z\|_{M_{\phi \psi }^*} = \sup_{t \in I} \phi (t) \psi (t)\,z^*(t)  \leq \sup_{t \in I} \phi (t) \psi (t)\,x^*(t) y^*(t) 
\leq \| x\|_{M_{\phi }^*}\, \| y\|_{M_{\psi }^*}.
$$
Therefore,
$$
\| z\|_{M_{\phi \psi }^*} \leq \inf \{ \| x\|_{M_{\phi }^*} \| y \|_{M_{\psi }^*}: z^* \leq x^* y^*\} \leq 2 \| z\|_{ M_{\phi }^* \odot M_{\psi }^*}.
$$

(ii) Let $z \in \Lambda_{\phi}  \odot M_{\psi }^*$. Then for any $\varepsilon > 0$ we can find 
$x \in \Lambda_{\phi},  y \in M_{\psi }^*$ such that $z = xy$ and $\| x \|_{\Lambda_{\phi}} \, \| y\|_{M_{\psi }^*} \leq (1+\varepsilon) \| z \|_{ \Lambda_{\phi}  \odot M_{\psi }^*}$.
Since $\psi^{\prime}(t) \leq \psi(t)/t$ and $\phi(t) \leq \frac{1}{a}\, t \phi^{\prime}(t)$ it follows that
\begin{eqnarray*}
\int_I z^*(t)\, d(\phi \psi)(t)
&\leq&
\int_I x^*(t/2) \, y^*(t/2) \left [\phi^{\prime}(t) \psi(t) + \phi(t) \psi^{\prime}(t) \right] dt\\
&\leq&
2\, \int_I x^*(t/2) \, y^*(t/2)\, \psi(t/2) \phi^{\prime}(t) \, dt \\
&+& 
\frac{1}{a}  \int_I x^*(t/2) \, y^*(t/2)\, t \phi^{\prime}(t) \frac{\psi(t)}{t} dt\\
&\leq&
2\, \sup_{s \in I} \psi(s) y^*(s)\, \int_I x^*(t/2) \, \phi^{\prime}(t) \, dt \\
&+& 
\frac{2}{a}\, \sup_{s \in I} \psi(s) y^*(s)\, \int_I x^*(t/2) \, \phi^{\prime}(t)\, dt\\
&\leq&
(2 + 2/a) \| y\|_{M_{\psi}^*}\, \| D_2 x\|_{\Lambda_{\phi}} \\
&\leq&
(4 + 4/a) \| y\|_{M_{\psi}^*}\, \| x\|_{\Lambda_{\phi}} \leq (4 + 4/a) (1+ \varepsilon) \, \| z \|_{  \Lambda_{\phi}  \odot M_{\psi }^*}.
\end{eqnarray*}
Since $\varepsilon > 0$ is arbitrary the first inclusion of (ii) is proved. To prove the second inclusion assume that 
$z = z^* \in  \Lambda_{\phi \psi }$. Then
$$
w(t) = z(t) \, \frac{\phi(t) \psi(t)}{t} \in L^1(I) ~ {\rm and} ~w = w^*.
$$
Moreover,
$$
\| w \|_{L^1} = \int_I z^*(t)  \frac{\phi(t) \psi(t)}{t} dt \leq \frac{1}{b} \, \int_I z^*(t) \, d(\phi \psi)(t) = \frac{1}{b} \|z \|_{\Lambda_{\phi \psi}}.
$$
Using the Lorentz result on the duality $(\Lambda_{\phi})^{\prime} \equiv M_{t/\phi(t)}$ (see [Lo51, Theorem 6], [Lo53, 
Theorem 3.6.1]; see also [KPS82, Theorem 5.2] for separable $\Lambda_{\phi}$, [HM92, Proposition 2.5(a)] with 
$p = q = 1$, [KM07, Theorem 2.2]) and Lozanovski{\u \i}'s factorization theorem 
$L^1 \equiv \Lambda_{\phi} \odot  ( \Lambda_{\phi})^{\prime}  \equiv  \Lambda_{\phi} \odot  M_{t/\phi(t)}$ 
we can find $u \in \Lambda_{\phi}, ~ v \in   M_{t/\phi(t)}$ such that
$$
w^* = u v ~~ {\rm and} ~~  \| u\|_{ \Lambda_{\phi}} \, \| v \|_{M_{t/\phi(t)}} \leq \| w\|_{L^1}.
$$
By Corollary 7 we can find $u_0 \in \Lambda_{\phi}, ~ v_0 \in M_{t/\phi(t)}  \overset{1}\hookrightarrow M_{t/\phi(t) }^*$ 
such that
$$
w^* \leq u_0^* v_0^* ~~ {\rm and} ~~  \| u_0 \|_{ \Lambda_{\phi}} \, \| v_0 \|_{M_{t/\phi(t)}} \leq 2 \,   \| u\|_{ \Lambda_{\phi}} \, \| v \|_{M_{t/\phi(t)}}.
$$
Let 
$$
x(t) = \frac{t}{\phi(t)} \frac{w^*(t)}{ \| v_0 \|_{M_{t/\phi(t)}}} ~~{\rm and} ~~ y = \frac{w^*}{u_0^*}.
$$
Then $x(t) \leq \frac{w^*(t)}{v_0^*(t)} \leq u_0^*(t)$ because $M_{t/\phi(t)}  \overset{1}\hookrightarrow M_{t/\phi(t) }^*$. 
Also $y(t) \leq v_0^*(t)$ and $z(t) \frac{\phi(t) \psi(t)}{t} = w(t) = x(t)\, \frac{\phi(t)}{t}  \| v_0 \|_{M_{t/\phi(t)}}$,
hence
$$
z(t) = x(t)\, \frac{1}{\psi(t)}\, \| v_0 \|_{M_{t/\phi(t)}}
$$
with $x \in \Lambda_{\phi}$ and $\frac{ \| v_0 \|_{M_{t/\phi(t)}}}{\psi(t)} \in M_{\psi}^*$.
Moreover,
\begin{eqnarray*}
\| z \|_{ \Lambda_{\phi}  \odot M_{\psi }^*} 
&\leq& 
\| x \|_{ \Lambda_{\phi}} \, \| \frac{1}{\psi} \|_{M_{\psi}^*} \, \| v_0\|_{M_{t/\phi(t)}} \leq
\| u_0 \|_{ \Lambda_{\phi}} \, \| v_0\|_{M_{t/\phi(t)}} \\
&\leq& 
2 \, \| u\|_{ \Lambda_{\phi}} \, \| v \|_{M_{t/\phi(t)}} \leq 2 \, \| w\|_{L^1} \leq \frac{2}{b} \, \| z \|_{\Lambda_{\phi \psi}}
\end{eqnarray*}
and the proof of (ii) is complete.

(iii) If $0 < p_{\phi, I}  \leq q_{\phi, I} < 1$, then both operators $H, H^*$ are bounded on 
$L^1(\frac{\phi(t)}{t})$ (see [KMP07], Theorem 4) and using Lemmas 4 and 5 we have
\begin{eqnarray*}
\Lambda_{\phi, 1}^{\theta} \, \Lambda_{\psi, 1}^{1-\theta} 
&=&
[L^1(\frac{\phi(t)}{t})^{(*)}]^{\theta}\, [L^1(\frac{\psi(t)}{t})^{(*)}]^{1-\theta} =
[L^1(\frac{\phi(t)}{t})^{\theta} \, [L^1(\frac{\psi(t)}{t})^{1-\theta}]^{(*)}\\
&=& 
L^1(\frac{\phi(t)^{\theta} \psi(t)^{1-\theta}}{t})^{(*)} = \Lambda_{\phi^{\theta} \psi^{1-\theta}, 1}
\end{eqnarray*}
with equivalent quasi-norms. Thus by Theorem 1(iv) we obtain
$$
\Lambda_{\phi, 1 } \odot \Lambda _{\psi, 1 } = (\Lambda_{\phi, 1 }^{1/2} \Lambda_{\psi, 1 }^{1/2})^{(1/2)}
= ( \Lambda_{\phi^{1/2} \psi^{1/2}, 1})^{(1/2)} = \Lambda _{\phi \psi, 1/2 }
$$
with equivalent quasi-norms. The last space is not normable since it contains isomorphic 
copy of $l^{1/2}$ (see [KM04], Theorem 1).

If $0 < p_{\phi, I}  \leq q_{\phi, I} < 1$, then both operators $H, H^*$ are bounded on $L^{\infty}(\phi)$ 
which can be proved directly. To show this we only need here to see equivalence of the corresponding integral 
inequalities on $\phi$ with assumptions on indices of $\phi$ and this is proved, for example, in [Ma85, 
Theorem 6.4] or [Ma89, Theorem 11.8] (see also [KPS], pp. 56-57). Then, using Lemmas 4 and 5, we have
\begin{eqnarray*}
\Lambda_{\phi, 1}^{\theta} \, (M_{\psi}^*)^{1-\theta} 
&=&
[L^1(\frac{\phi(t)}{t})^{(*)}]^{\theta}\, [L^{\infty} (\psi)^{(*)}]^{1-\theta} =
[L^1(\frac{\phi(t)}{t})^{\theta} \, [L^{\infty}(\psi)^{1-\theta}]^{(*)}\\
&=& 
L^{1/{\theta}} (\frac{\phi(t)^{\theta} \psi(t)^{1-\theta}}{t^{\theta}})^{(*)} = \Lambda_{\phi^{\theta} \psi^{1-\theta}, 1/{\theta}}
\end{eqnarray*}
with equivalent quasi-norms. Thus, by Theorem 1(iv), we obtain
$$
\Lambda_{\phi, 1 } \odot M_{\psi }^* = [\Lambda_{\phi, 1 }^{1/2} (M_{\psi}^*)^{1/2}]^{(1/2)} =
( \Lambda_{\phi^{1/2} \psi^{1/2}, 2})^{(1/2)} = \Lambda _{\phi \psi, 1 }
$$
with equivalent quasi-norms. Similarly for Marcinkiewicz spaces
\begin{eqnarray*}
(M_{\phi}^*)^{\theta} \, (M_{\psi}^*)^{1-\theta} 
&=&
[L^{\infty}(\phi)^{(*)}]^{\theta}\, [L^{\infty} (\psi)^{(*)}]^{1-\theta} =
[L^{\infty}(\phi)^{\theta} \, L^{\infty}(\psi)^{1-\theta}]^{(*)}\\
&=& 
L^{\infty} (\phi^{\theta} \psi^{1-\theta})^{(*)} = M_{\phi^{\theta} \psi^{1-\theta}}^*
\end{eqnarray*}
and, by Theorem 1(iv), we obtain
$$
M_{\phi}^* \odot M_{\psi }^* = [(M_{\phi}^*)^{1/2} (M_{\psi}^*)^{1/2}]^{(1/2)} =
( M_{\phi^{1/2} \psi^{1/2}}^*)^{(1/2)} = M_{\phi \psi }^*
$$
with equivalent quasi-norms. This proves theorem completely.
\end{proof}

\begin{center}
\textbf{6. Factorization of some Banach ideal spaces}
\end{center}

The factorization theorem of Lozanovski\u{\i} states that for any Banach ideal space $E$ 
the space $L^1$ has a factorization $L^1 \equiv E \odot E^{\prime }$. The natural generalization 
of the type 
\begin{equation} \label{factorization}
F \equiv E \odot M(E, F)
\end{equation}
is not true without additional assumptions on the spaces, as we can see on the example below.

\vspace{3mm} 
{\bf Example 2}. If $E = L^{p, 1}$ with the norm $\| x\|_E = \frac{1}{p} \int_I t^{1/p-1} x^*(t) \, dt$ 
for $1 < p < \infty$, then $M(L^{p, 1}, L^p) \equiv L^{\infty}$ (cf. [MP89], Theorem 3) and
$$
L^{p, 1} \odot M(L^{p, 1}, L^p) \equiv L^{p, 1} \odot L^{\infty} \equiv L^{p, 1} \subsetneq L^p.
$$
Therefore, factorization (\ref{factorization}) is not true and we even don't have factorization 
$L^p = E \odot M(E, L^p)$ with equivalent norms. Similarly, if $F = L^{p, \infty}$ with the norm 
$\| x\|_F = \sup_{t \in I} t^{1/p} x^{**}(t)$ for $1 < p < \infty$, then 
$M(L^p, L^{p, \infty}) \equiv M(L^{p^{\prime}, 1}, L^{p^{\prime}}) \equiv L^{\infty}$ and
$$
L^p \odot M(L^p, L^{p, \infty}) \equiv L^p \odot L^{\infty} \equiv L^p \subsetneq L^{p, \infty}.
$$
Therefore, again factorization (\ref{factorization}) is not true and we even don't have equality 
$F = L^p \odot M(L^p, F)$ with equivalent norms.
\vspace{2mm}

Let us collect some factorization results of type (\ref{factorization}). First of all the Lozanovski{\u \i} 
factorization theorem was announced in 1967 (cf. [Lo67], Theorem 4) and published with detailed 
proof in 1969 (cf. [Lo69], Theorem 6). His proof uses the Calder\'on space $F = E^{1/2} ( E^{\prime })^{1/2}$ 
and result about its dual $F^{\prime \prime} \equiv F^{\prime} \equiv L^2$ (cf. [Lo69], Theorem 5; 
see also [Ma89, p. 185] and [Re93]). 
Lozanovski{\u \i}  factorization theorem was new even for finite dimensional spaces. 
In 1976 Jamison and Ruckle [JR76] proved that $l^1$ factors through every normal Banach 
sequence space and its K\"othe dual. Proof even in the finite dimensional case is indirect and it 
uses the Brouwer fixed point theorem.
Later on Lozanovski{\u \i}'s factorization result was proved by Gillespie [Gi81], using different method, 
which inspiration was coming from the theory of reflexive algebras of operators on 
Hilbert space.

If $E, F$ are finite dimensional ideal spaces and $B_E, B_F$ denote their unit balls, then 
Bollob\'as and Leader [BL95], with the help of Jamison-Ruckle method, proved factorization 
$B_E  \odot B_{M(E, F)} \equiv B_F$ under assumptions that $B_F$ is a strictly unconditional 
body and $B_{M(E, F)} $ is smooth. 
\vspace{1mm}

Nilsson [Ni85, Lemma 2.5], using the Maurey factorization theorem (cf. [Ma74, Theorem 8]; see also 
[Wo91, pp. 264-266]), proved the following result of type (\ref{factorization}): 
if $E$ is a Banach ideal space which is $p$-convex with constant $1$, 
then 
\begin{equation} \label{equation29}
E^{\prime} \equiv L^{p^{\prime}} \odot M(E, L^p)  \equiv L^{p^{\prime}} \odot M(L^{p^{\prime}}, E^{\prime}).
\end{equation}
By duality result and (\ref{equation29}) we obtain that if $F$ is a Banach ideal space with the Fatou property 
which is $q$-concave with constant $1$ for $ 1 < q < \infty$, then 
\begin{equation} \label{equation30}
F = F^{\prime \prime} \equiv L^q \odot M(F^{\prime}, L^{q^{\prime}}) \equiv L^q \odot M(L^q, F).
\end{equation}
Factorization (\ref{equation30}) was proved and used by Nilsson [Ni85, Theorem 2.4] in a new proof of 
the Pisier theorem (cf. [Pi79a, Theorem 2.10], [Pi79b, Theorem 2.2]; see also [TJ89, Theorem 28.1]): if a Banach 
ideal space $E$ with the Fatou property is $p$-convex and $q$-concave with constants $1, 1 < r < \infty$, 
and $\frac{1}{p} = \frac{\theta}{r} + 1-\theta, \frac{1}{q} = \frac{\theta}{r}$, then the space 
$E_0 \equiv M(L^q, E)^{(1-\theta)}$ is a Banach ideal space and $E \equiv E_0^{1-\theta} (L^r)^{\theta}$.
First part of the proof follows from the facts that
\begin{equation*}
\| |x|^{1-\theta} \|_{M(L^q, E)}^{\frac{1}{1-\theta}} \equiv \| |x|^{1/s^{\prime}} \|_{M(L^s, E^{(1/p)})}^{s^{\prime}} 
\equiv \| |x|^{1/s^{\prime}} \|_{M( (E^{(1/p)})^{\prime}, L^{s^{\prime}})}^{s^{\prime}},
\end{equation*}
$E^{(1/p)}$ is a Banach space and $M( (E^{(1/p)})^{\prime}, L^{s^{\prime}})$ is $s^{\prime}$-convex with constant 
$1$, where $s = \frac{r}{\theta p}$. Second part uses (\ref{equation30}) and by Theorem 1(ii) we obtain
\begin{equation*}
E_0^{1-\theta} (L^r)^{\theta} \equiv E_0^{\frac{1}{1-\theta}} \odot L^{\frac{r}{\theta}} \equiv E_0^{\frac{1}{1-\theta}} \odot L^q \equiv M(L^q, E) \odot L^q \equiv E.
\end{equation*}

Schep proved factorization (\ref{equation30}) and also the reverse implication, that is, if (\ref{equation30}) holds, then the space 
$F$ is $q$-concave with constant $1$ (cf. [Sc10], Theorem 3.9). He has also proved another factorization result 
(even equivalence - see [Sc10], Theorem 3.3): if Banach ideal space $E$ with the Fatou property 
is $p$-convex with constant $1 \, {\rm (}1 < p < \infty {\rm )}$, then
\begin{equation} \label{equation31}
L^p \equiv E \odot M(E, L^p).
\end{equation}
His proof has misprints in Theorem 3.2. The proof should be as follows: using property (g) from [MP89] we obtain
\begin{equation*}
M(E, L^p)^{(1/p)} \equiv M(E^{(1/p)}, (L^p)^{(1/p)}) \equiv M(E^{(1/p)}, L^1) \equiv [E^{(1/p)}]^{\prime},
\end{equation*}
and by the Lozanovski{\u \i} factorization theorem
\begin{equation*}
E^{(1/p)} \odot M(E, L^p)^{(1/p)} \equiv E^{(1/p)} \odot [E^{(1/p)}]^{\prime} \equiv L^1. 
\end{equation*}
Taking then $p$-convexification on both sides and using Theorem 1(iii) we get
\begin{equation*}
E \odot M(E, L^p) \equiv [E^{(1/p)} \odot M(E, L^p)^{(1/p)}]^{(p)} \equiv  (L^1)^{(p)} \equiv L^p.
\end{equation*}

Note that factorization theorem of the type (\ref{equation30}): $F = l^q \odot M(l^q, F)$ for any $q$-concave 
Banach space $F$ with a monotone unconditional basis was proved already in 1980 (cf. [LT-J80], Corollary 3.2).
\vspace{2mm}

If a space $E$ has the Fatou property, then in the definition of the norm of $E \odot E^{\prime }$ 
we may take ``minimum" instead of ``infimum". It is known that the Fatou property of $E$ 
is equivalent with the isometric equality $E \equiv E^{\prime \prime}$. Then $E$ is called 
{\it perfect}. This notion can be generalized to $F$-perfectness. We say that $E$ is 
{\it $F$-perfect} if $M(M(E, F), F) \equiv E$ (see [MP89], [CDS08] and [Sc10] for more information 
about $F$-perfectness). Is there any connection between factorization (\ref{factorization}) 
and to be $F$-perfect by $E$?

\vspace{3mm}
{\bf Theorem 8.} {\it Let $E, F$ be Banach ideal spaces with the Fatou property. Then 
factorization $E\odot M(E, F) \equiv F$ implies $F$-perfectness of $E$, i.e., $M( M( E, F) ,F) \equiv E$.}

\begin{proof} Schep [Sc10, Theorem 2.8] proved that if $E\odot F$ is a Banach ideal space, then
$M(E, E \odot F) \equiv F$ (see also Theorem 4 above). Since $E \odot M(E, F) \equiv F$ is a Banach ideal space by assumption,
therefore from the above Schep result we obtain
\begin{equation*}
M( M(E, F), F) \equiv M( M(E, F), E \odot M(E, F)) \equiv E,
\end{equation*} 
\vspace{-2mm}
which is $F$-perfectness of $E$.
 \end{proof}
The example of Bollobas and Brightwell [BB00], presented in [Sc10, Example 3.6], shows 
that the reverse implication is not true, even for three-dimensional spaces.
\vspace{2mm}

Almost all proofs in factorization theorems are tricky or use powerful theorems and, in fact,
equality $E \odot M(E, F) \equiv F$ is proved without calculating $M(E, F)$
directly. Except some special cases it seems to be the only way to prove
equality of the norms in (\ref{factorization}). However, it seems to be also
useful to have equality (\ref{factorization}) with just equivalence of the norms, 
that is, 
\begin{equation} \label{equation32}
F = E \odot M(E, F).
\end{equation}
This can be done by finding $M(E, F)$ and $E \odot M(E, F)$ separately 
and we will do so. Observe also that if a Banach ideal space $E$ is $p$-convex ($1 < p < \infty$) 
with constant $K > 1$, then $E^{(1/p)}$ is $1$-convex with constant $K^p$ and
$$
\| x \|^0 = \inf \{\sum_{k=1}^n \| x_k \|_{E^{(1/p)}}: |x | \leq \sum_{k=1}^n |x_k |, x_k \in E^{(1/p)}, n \in {\mathbb N} \}
$$
defines norm on $E^{(1/p)}$ with $K^{-p} \| x \|_{E^{(1/p)}} \leq \| x \|^0 \leq  \| x \|_{E^{(1/p)}}$. Thus $E_0 = (E^{(1/p)}, \| \cdot \|^0)$ is a Banach ideal space and its $p$-convexification $E_0^{(p)} = E$ with the norm $\| x\|^1 = [\| |x|^p \|^0)^{1/p}$ is $p$-convex with constant $1$ (cf. [LT79, Lemma 1.f.11] and [ORS08, Proposition 2.23]), and we can use result from (\ref{equation29}) to obtain $E^{\prime} = L^{p^{\prime}} \odot M(L^{p^{\prime}}, E^{\prime})$.
\vspace{1mm}

As a straightforward conclusion from Corollary 6.1 in [KLM12] and Theorem A(a)
with Theorem 5(a) we get the following factorization theorem for the Calder\'{o}n-Lozanovski{\u \i}  
$E_{\varphi}$-spaces.

\vspace{3mm}
{\bf Theorem 9.} {\it Let $E$ be a Banach ideal space with the Fatou property and ${\it supp} E = \Omega$. 
Suppose that for two Young functions $\varphi ,\varphi_1$ there is a Young function $\varphi _{2}$
such that one of the following conditions holds:
\begin{itemize}
\item[$(i)$]  $\varphi _{1}^{-1}\varphi _{2}^{-1}\approx \varphi ^{-1}$ for all
arguments,
\item[$(ii)$] $\varphi _{1}^{-1}\varphi _{2}^{-1}\approx \varphi ^{-1}$ for large
arguments and $L^{\infty }\hookrightarrow E$,
\item[$(iii)$] $\varphi _{1}^{-1}\varphi _{2}^{-1}\approx \varphi ^{-1}$ for small
arguments and $E\hookrightarrow L^{\infty }$.
\end{itemize}
Then the factorization $E_{\varphi _{1}}\odot M\left( E_{\varphi _{1}},E_{\varphi }\right)
=E_{\varphi }$ with equivalent norms is valid and, in consequence, the space $E_{\varphi _{1}}$ is 
$E_{\varphi }$--perfect up to equivalence of norms.}
\vspace{2mm}

Moreover, applying Lemma 7.4 from [KLM12] to the Theorem 9(i) one has the following special case.

\vspace{3mm}
{\bf Corollary 8.} {\it Let $\varphi ,\varphi _{1}$ be two Orlicz functions, and let $E$ Banach an ideal 
space with the Fatou property and ${\it supp} E=\Omega$. If the function
$f_{v}(u) = \frac{\varphi (uv)}{\varphi _1(u)}$ is non-increasing on $(0,\infty )$ for any $v>0$, then 
the factorization $E_{\varphi _1} \odot M( E_{\varphi _1},E_{\varphi }) = E_{\varphi }$ 
is valid with equivalent norms and, in consequence, the space $E_{\varphi _1}$ is $E_{\varphi }$--perfect 
up to equivalence of the norms.}

\begin{proof} It is enough to take as $\varphi_2$ the function defined by
\vspace{-2mm}
$$
\varphi _2(u) = (\varphi \ominus \varphi _1)(u) = \sup_{v>0}[\varphi (uv)-\varphi _1(v)]
$$ 
and use the fact proved in [KLM12, Lemma 7.4] showing that 
$\varphi _{1}^{-1}\varphi _{2}^{-1}\approx \varphi ^{-1}$ for all arguments.
\end{proof}

Before we consider factorization of Lorentz and Marcinkiewicz let us calculate 
``missing" multipliers spaces.

\vspace{3mm}
{\bf Proposition 3.} {\it Suppose $\phi , \psi $ are non-decreasing, concave functions 
on $I$ with $\phi ( 0^+) = \psi (0^+) =0$. Let $E$ and $F$ be symmetric spaces on $I$ with 
fundamental functions $f_E(t) = \phi(t)$ and $f_F(t) = \psi(t)$. If 
$\omega ( t) = \sup_{0 < s \leq t} \frac{\psi (s) }{\phi ( s) }$ is finite for any $t \in I$, then
\vspace{-2mm}
\begin{equation*}
M(\Lambda_{\phi}, F) \overset{1}{\hookrightarrow } M_{\omega}, ~~
M(E, M_{\psi}) \overset{1}{\hookrightarrow } M_{\omega},~~ {\rm and} ~~
M_{\omega}^* \overset{1}{\hookrightarrow } M(\Lambda_{\phi, 1}, \Lambda_{\psi, 1}).
\end{equation*}

\vspace{-3mm}

\noindent
If, moreover, $s_{\phi, I} \geq a > 0$, then $M_{\omega}^* \overset{1/a}{\hookrightarrow } 
M(\Lambda_{\phi}, \Lambda_{\psi})$.}

\begin{proof} By Theorem 2.2(iv) in [KLM12] the function $\omega$ is a fundamental function of 
$M(\Lambda_{\phi}, F)$ and by the maximality of the space $M_{\omega }$ and Theorem 2.2(i) in [KLM12] 
we obtain imbedding $M(\Lambda_{\phi}, F) \overset{1}{\hookrightarrow } M_{\omega} $. 
On the other hand, using the property (e) from [MP89, p. 326] about duality of multipliers, the duality 
$(M_{\psi})^{\prime} \equiv \Lambda_{t/{\psi(t)}}$ and the above result we obtain
\vspace{-2mm}
\begin{equation*}
M(E, M_{\psi}) \equiv M((M_{\psi})^{\prime}, E^{\prime}) \equiv M(  \Lambda_{t/{\psi(t)}}, 
E^{\prime}) \overset{1}{\hookrightarrow } M_{\omega}.
 \end{equation*}
Two other imbeddings will be proved if we show that $\frac{1}{\omega}$ belongs to 
the corresponding spaces. Since, by Theorem 2.2(ii) in [KLM12], we have 
$\| y \|_{M(E, F)} = \sup_{\| x^*\|_E \leq 1} \| x^* y^* \|_F$ it follows that
\begin{eqnarray*}
\| \frac{1}{\omega}\|_{M(\Lambda_{\phi, 1}, \Lambda_{\psi,1})} 
&=& 
\sup_{\| x \|_{\Lambda_{\phi, 1}} \leq 1} \int_I \left (x^* \, \frac{1}{\omega} \right)^*(t) \frac{\psi(t)}{t} \, dt 
=
\sup_{\| x \|_{\Lambda_{\phi, 1}} \leq 1} \int_I x^*(t) \frac{1}{\omega(t)} \frac{\psi(t)}{t} \, dt\\
&=&
\sup_{\| x \|_{\Lambda_{\phi, 1}} \leq 1} \int_I x^*(t) \inf_{0 < s \leq t} \frac{\phi(s)}{\psi(s)} \frac{\psi(t)}{t} \, dt 
\leq \sup_{\| x \|_{\Lambda_{\phi, 1}} \leq 1} \int_I x^*(t) \frac{\phi(t)}{t} \, dt \leq 1,
\end{eqnarray*}
\vspace{-2mm}
and, again by the above mentioned result in [KLM12],
\begin{eqnarray*}
\| \frac{1}{\omega}\|_{M(\Lambda_{\phi}, \Lambda_{\psi})} 
&=& 
\sup_{\| x \|_{\Lambda_{\phi}} \leq 1} \int_I x^*(t) \frac{1}{\omega(t)} \psi^{\prime}(t) \, dt 
= \sup_{\| x \|_{\Lambda_{\phi}} \leq 1} \int_I x^*(t) \inf_{0 < s \leq t} \frac{\phi(s)}{\psi(s)} \, \psi^{\prime}(t) \, dt\\
&\leq&
\sup_{\| x \|_{\Lambda_{\phi}} \leq 1} \int_I x^*(t) \frac{\phi(t)}{t} \, dt 
\leq \sup_{\| x \|_{\Lambda_{\phi}} \leq 1} \frac{1}{a} \, \int_I x^*(t) \,\phi^{\prime}(t)\, dt \leq 1/a,
\end{eqnarray*}
and all imbeddings are proved. \end{proof}

Putting together previous results on products and multipliers of Lorentz and Marcinkie\-wicz spaces 
we are ready to proof factorization of these spaces.

\vspace{2mm}
{\bf Theorem 10.} {\it Let $\phi ,\psi $ be a non-decreasing, concave functions on $I$ 
with $\phi (0^+) = \psi (0^+) = 0$. Suppose $\frac{\psi(t)}{\phi(t)}$ is a non-decreasing function on $I$.
\vspace{-2mm}
\begin{itemize}
\item[$(a)$] If $s_{\phi, I} > 0$ and $s_{\psi, I} > 0$, then $\Lambda_{\phi} \odot M(\Lambda_{\phi}, \Lambda_{\psi}) 
= \Lambda_{\psi}$. 
\end{itemize}
\vspace{-1mm}
Moreover, for any symmetric space $F$ on $I$ with the fundamental function $f_F(t) = \psi(t)$ and under the 
above assumptions on $\phi$ and $\psi$ we have
\begin{equation} \label{equation32}
\Lambda_{\phi} \odot M(\Lambda_{\phi}, F) = F ~ {\it if ~and ~only ~if} ~~ F = \Lambda_{\psi}.
\end{equation}
\begin{itemize}
\item[$(b)$] If $\sigma_{\phi, I} < 1$ and $\sigma_{\psi, I} < 1$, then $M_{\phi} \odot M(M_{\phi}, M_{\psi}) = M_{\psi}$. 
\end{itemize}
\vspace{-1mm}
Moreover, for any symmetric space $E$ on $I$ having Fatou property, with the fundamental function $f_E(t) = \phi(t)$ and under the above assumptions on $\phi$ and $\psi$ we have
\begin{equation}  \label{equation33}
E \odot M(E, M_{\psi}) = M_{\psi} ~ {\it if ~and ~only ~if} ~~ E =M_{\phi}.
\end{equation}

\vspace{-5mm}

\begin{itemize}
\item[$(c)$] If $\sigma_{\phi, I} < 1, s_{\psi, I} > 0$ and $s_{\psi/\phi, I} > 0$, then 
\begin{equation}  \label{equation34}
M_{\phi} \odot M(M_{\phi}, \Lambda_{\psi}) = \Lambda_{\psi}.
\end{equation} 
\end{itemize}
}
\vspace{-3mm}

\begin{proof} (a) Using Proposition 3 we have
\begin{equation*}
M(\Lambda_{\phi}, \Lambda_{\psi} ) \overset{1}{\hookrightarrow } M_{\omega} \overset{1}{\hookrightarrow } M_{\omega}^* \overset{1/a}{\hookrightarrow } M(\Lambda_{\phi}, \Lambda_{\psi}), ~ {\rm where} ~ \omega(t) = \sup_{0 < s \leq t} \frac{\psi(s)}{\phi(s)}.
\end{equation*}
Since $\psi/\phi$ is a non-decreasing function on $I$ it follows that $\phi \omega = \psi, s_{\phi \omega; I} = s_{\psi, I} > 0$ and, by Theorem 7(ii),
\begin{equation*}
\Lambda_{\phi} \odot M(\Lambda_{\phi}, \Lambda_{\psi}) = \Lambda_{\phi} \odot M_{\omega}^* = \Lambda_{\psi}
\end{equation*}
with equivalent norms. Under the assumptions on $F$ we have from Proposition 3 the imbedding $M(\Lambda_{\phi}, F) \overset{1}{\hookrightarrow } M_{\omega}$ and then, by Theorem 7(ii),
\begin{equation*}
F = \Lambda_{\phi} \odot M(\Lambda_{\phi}, F) \overset{1}{\hookrightarrow } \Lambda_{\phi} \odot M_{\omega} \overset{1}{\hookrightarrow } \Lambda_{\phi} \odot M_{\omega}^* = \Lambda_{\psi}.
\end{equation*}
Minimality of $\Lambda_{\psi}$ gives $F = \Lambda_{\psi}$.
\vspace{2mm}

(b) From the fact that $(M_{\phi})^{\prime} \equiv \Lambda_{t/\phi(t)}$, the general duality property of multipliers 
(see [MP89], property (e)) and using Proposition 3 we obtain
\begin{equation*}
M(M_{\phi}, M_{\psi}) \equiv M(M_{\psi}^{\prime}, M_{\phi}^{\prime}) \equiv M(\Lambda_{t/\psi(t)}, \Lambda_{t/\phi(t)}) = M_{\omega}^*
\end{equation*}
because $s_{t/\psi(t), I} = 1 - \sigma_{\psi, I} > 0$. Since $\psi/\phi$ is a non-decreasing function on $I$ it follows that $\phi \omega = \psi, \sigma_{\phi \omega; I} = \sigma_{\psi, I} < 1$ and by Theorem 7(i) with the fact that $\sigma_{\phi, I} < 1$ we have
\begin{equation*}
M_{\phi} \odot M(M_{\phi}, M_{\psi}) = M_{\phi} \odot M_{\omega} = M_{\phi}^* \odot M_{\omega}^* 
= M_{\phi \omega}^* = M_{\psi}^* = M_{\psi}.
\end{equation*}
Under the assumptions on $E$ we obtain from Proposition 3 that $ M(E, M_{\psi}) \overset{1}{\hookrightarrow } 
M_{\omega}$ and 
\begin{equation*}
M_{\psi} = E \odot M(E, M_{\psi}) \overset{1}{\hookrightarrow }  E \odot M_{\omega}.
\end{equation*}
On the other hand, by Theorem 7(i) and assumption $q_{\psi; I} < 1$
\begin{equation*}
M_{\phi} \odot M_{\omega}  \overset{1}{\hookrightarrow } M_{\phi}^* \odot M_{\omega}^* \overset{2}{\hookrightarrow } M_{\phi \omega}^* \equiv M_{\psi}^* = M_{\psi}.
\end{equation*}
Therefore, $M_{\phi} \odot M_{\omega}  \overset{C}{\hookrightarrow } E \odot M_{\omega}$. Using now Schep's theorem, saying that if $E \odot F  \overset{C}{\hookrightarrow } E \odot G$, then $F  \overset{C}{\hookrightarrow } G$ (see [Sc10], Theorem 2.5), we obtain $M_{\phi}  \overset{C}{\hookrightarrow } E$. Maximality of Marcinkiewicz space $M_{\phi}$ implies that $E = M_{\phi}$ since fundamental function of $E$ is $f_E(t) = \phi(t)$ for all $t \in I$.

(c) Using Theorem 2.2(v) from [KLM12] we obtain $M(M_{\phi}, \Lambda_{\psi}) \equiv \Lambda_{\eta}$, where $\eta(t) = \int_0^t (\frac{s}{\phi(s)} )^{\prime} \psi^{\prime}(s)\, ds < \infty$. Since
\begin{eqnarray*}
\eta(t) 
&=&
\int_0^t \frac{\phi(s) \psi^{\prime}(s) - \phi^{\prime}(s) s \psi^{\prime}(s)}{\phi(s)^2} \, ds 
\leq \int_0^t \frac{ \psi^{\prime}(s)}{\phi(s)} \, ds\\
&\leq&
\frac{1}{ s_{\psi/\phi}} \, \int_0^t (\frac{\psi}{\phi})^{\prime}(s) \, ds = \frac{1}{ s_{\psi/\phi}} \, \frac{\psi(t)}{\phi(t)},
\end{eqnarray*}
and
\begin{eqnarray*}
\eta(t) 
&=&
\int_0^t \frac{\phi(s) \psi^{\prime}(s) - \phi^{\prime}(s) s \psi^{\prime}(s)}{\phi(s)^2} \, ds \geq
\int_0^t \frac{\phi(s) \psi^{\prime}(s) - \phi^{\prime}(s) \psi(s)}{\phi(s)^2} \, ds\\
&=&
 \int_0^t (\frac{\psi}{\phi})^{\prime}(s) \, ds = \frac{\psi(t)}{\phi(t)},
\end{eqnarray*}
it follows that $M(M_{\phi}, \Lambda_{\psi})= \Lambda_{\psi/\phi}$. Using to this equality assumption 
$\sigma_{\phi, I} < 1$ and result from Theorem 7(ii) we obtain
\begin{equation*}
M_{\phi} \odot M(M_{\phi}, \Lambda_{\psi}) = M_{\phi} \odot  \Lambda_{\psi/\phi} = M_{\phi}^* \odot  \Lambda_{\psi/\phi} = \Lambda_{\psi},
\end{equation*} 
and the theorem is proved. \end{proof}
\vspace{2mm}

Applying the above theorem to classical Lorentz $L^{p, 1}$ and Marcinkiewicz $L^{p, \infty}$ spaces we obtain 
the following factorization results:

\vspace{3mm} 
{\bf Example 3}. (a) If $1 \leq p \leq q < \infty$, then $L^{p, 1} = L^{q, 1} \odot M(L^{q, 1}, L^{p, 1})$. 

(b) If $1 < p \leq q \leq \infty$, then $L^{p, \infty} = L^{q, \infty} \odot M(L^{q, \infty}, L^{p, \infty})$. 

(c)  If $1 < p < q \leq \infty$, then $L^{p, 1} = L^{q, \infty} \odot M(L^{q, \infty}, L^{p, 1})$.
\vspace{2mm}

\noindent
What about factorization in classical Lorentz $L^{p, q}$-spaces?  

\vspace{3mm} 
{\bf Example 4}. If either $1 \leq r \leq p < q < \infty$ or $1 < p < q \leq r \leq \infty$, then 
\begin{equation*}
L^{p, r} = L^{q, r} \odot M(L^{q, r}, L^{p, r}).
\end{equation*} 
In fact, if  $1 \leq r \leq p < q < \infty$ then using the commutativity of $r$-convexification with 
multipliers (see property (g) in [MP89]) and Proposition 3 we obtain
\begin{equation*}
M(L^{q, r}, L^{p, r}) \equiv M((L^{q/r, 1})^{(r)}, (L^{p/r, 1})^{(r)}) \equiv M(L^{q/r, 1}, L^{p/r, 1})^{(r)} 
= (L^{pq/[r(q-p)], \infty})^{(r)}.
\end{equation*}
Finally, by Theorem 1(iii) and Theorem 7(ii) with $\phi(t) = t^{r/q}$ and $\psi(t) = t^{r/p-r/q}$, we obtain
\begin{eqnarray*}
L^{q, r} \odot M(L^{q, r}, L^{p, r}) 
&=& 
(L^{q/r, 1})^{(r)} \odot  (L^{pq/[r(q-p)], \infty})^{(r)} \\
&=&
(L^{q/r, 1} \odot  L^{pq/[r(q-p)], \infty})^{(r)} = (L^{p/r, 1})^{(r)} = L^{p, r}.
\end{eqnarray*}
The case $1 < p < q \leq r \leq \infty$ can be proved by duality of multipliers and the above calculations.

\vspace{3mm}
{\bf Theorem 11.} {\it Let $\phi $ be an increasing, concave function on $I$ with 
$0 < p_{\phi, I}  \leq q_{\phi, I} < 1$. 

\vspace{-2mm}

\begin{itemize}
\item[$(a)$] Suppose that $F$ is a symmetric space on $I$ with the lower Boyd index 
$\alpha_F > q_{\phi, I}$ and such that $M(M_{\phi}^*, F) \neq \{0\}$. Then 

\vspace{-2mm}

\begin{equation*}
F = M_{\phi}^* \odot M(M_{\phi}^*, F) = M_{\phi} \odot M(M_{\phi}, F).
\end{equation*}

\vspace{-2mm}

\item[$(b)$] Suppose that $E$ is a symmetric space on $I$ with the Fatou property, which Boyd 
indices satisfy $0 < \alpha_E \leq \beta_E < p_{\phi, I}$ and such that $M(E, \Lambda_{\phi}) \neq \{0\}$. 
Then 
\vspace{-2mm}
\begin{equation*}
\Lambda_{\phi, 1} = E \odot M(E, \Lambda_{\phi}).
\end{equation*}
\end{itemize}

}

\noindent
Let us start with the following identifications.

\vspace{3mm}
{\bf Lemma 6.} \label{Lemma6} {\it Under assumptions on $\phi$ from Theorem 11 we have
\begin{equation} \label{identification31}
M(L^{\infty}(\phi)^{(*)}, F) =  M(L^{\infty}(\phi), F)^{(*)} \equiv F(1/ \phi)^{(*)}.
\end{equation}

}

\begin{proof}
Since we have equivalences 
$$
z \in M(L^{\infty}(\phi), F)^{(*)} \Leftrightarrow z^* \in M(L^{\infty}(\phi), F)\Leftrightarrow \frac{z^*}{\phi} \in F \Leftrightarrow z^* \in F(1/\phi) \Leftrightarrow z \in F(1/ \phi)^{(*)}
$$ 
with equalities of the norms, the equality $M(L^{\infty }(\phi ), F)^{(\ast )}\equiv F(1/\phi )^{(\ast )}$ follows. Let 
$I = ( 0, 1)$. We prove the equality $M(L^{\infty }(\phi )^{(\ast )}, F) = F(1/\phi )^{(\ast )}.$
Clearly, it is enough to show that the following conditions are equivalent:
\vspace{1mm}

$1^0$ $z\in M(L^{\infty }(\phi )^{(\ast )}, F)$, 
\vspace{1mm}

$2^0$ $\frac{z}{\phi \left( \omega \right) }\in F$ for every measure preserving transformation 
$\left( mpt\right) $ $\omega :I\rightarrow I$, 
\vspace{1mm}

$3^0$ $\frac{z^{\ast }}{\phi }\in F$.

\noindent
Moreover, we prove that 
\begin{equation}  \label{equation36}
\| \frac{z^*}{\phi} \|_F \leq \sup_{\omega-{\rm mpt}} \| \frac{z}{\phi(\omega)} \|_F = \| z \|_{M[L^{\infty}(\phi)^{(*)}, F]} \leq \| D_2\|_{F \rightarrow F}\, \| \frac{z^*}{\phi} \|_F.
\end{equation}

$1^0 \Rightarrow 2^0$. Let $z \in M(L^{\infty}(\phi )^*, F)$ and take arbitrary mpt $\omega: I\rightarrow I.$ 
Since $\left( \frac{1}{\phi ( \omega) }\right)^* \phi = \frac{1}{\phi }\phi =1$ it follows that 
$\frac{1}{\phi ( \omega) }\in L^{\infty }(\phi )^{(\ast )}$. Whence 
$\frac{z}{\phi (\omega ) }\in F$ and $\| \frac{z}{\phi ( \omega ) } \| _{F} \leq \| z \|_{M[L^{\infty }(\phi )^{(*)}, F]}.$ 
Consequently,
$$
\sup_{\omega-{\rm mpt}} \left \| \frac{z}{\phi(\omega)} \right \|_F \leq  \| z \|_{M(L^{\infty }(\phi )^{(*)}, F)}.
$$

$2^0 \Rightarrow 1^0$. Let $x\in L^{\infty}(\phi )^{(* )} $ with the norm $\leq 1$. Take a mpt $\omega _{0}$ such
that $| x | = x^*( \omega _0)$. Then 
$$
| z x | = | z | \, x^*( \omega_0) \leq \frac{| z |}{\phi ( \omega_0) }\in F,
$$
because $1 \geq \| x^* \phi \| _{L^{\infty }} = \| ( x^* \phi) ( \omega _0) \|_{L^{\infty }}.$ 
Thus $z\in M(L^{\infty }(\phi )^{(* )}, F)$ and
$$
\| z \| _{M(L^{\infty }(\phi )^{(* )}, F)}  \leq  \sup\limits_{\omega -mpt} \| \frac{z}{\phi ( \omega) } \|_{F}.
$$

$2^0 \Rightarrow 3^0 .$ Take a mpt $\omega _0$ such that $| z | = z^* ( \omega_0) .$ Then
$$
\frac{z^*}{\phi }\sim \frac{z^*}{\phi } ( \omega _0) = \frac{| z |}{\phi ( \omega _0) }\in F \text{
and } \| \frac{z^*}{\phi } \| _{F} = \| \frac{ | z | }{\phi ( \omega _0) } \|_{F} 
\leq \sup\limits_{\omega -mpt} \| \frac{z}{\phi ( \omega) } \|_{F}.
$$

$3^0 \Rightarrow 2^0.$ For each mpt $\omega: I \rightarrow I$ we have 
$$
\frac{z}{\phi ( \omega ) }( t) \sim ( \frac{z}{\phi ( \omega ) } ) ^*( t) \leq 
z^*( t/2) ( \frac{1}{\phi ( \omega ) } )^*( t/2) = \frac{z^*( t/2) }{\phi(t/2) } = 
D_{2} ( \frac{z^*}{\phi })(t).
$$
By symmetry of $F$ we obtain $\frac{z}{\phi ( \omega ) }\in F$ and 
$$
\sup\limits_{\omega -mpt} \| \frac{z}{\phi ( \omega) } \|_{F} \leq 
\| D_{2} \| _{F\rightarrow F} \| \frac{z^*}{\phi } \|_{F}.
$$
The proof of (\ref{equation36}) and also (\ref{identification31}) is finished for $I = (0, 1)$. If $I = (0, \infty)$, then 
the existence of a measure preserving transformation $\omega_0: I \rightarrow I$ requires additional assumption, 
in the first case that $\phi(\infty) = \infty$, which we have because $p_{\phi, I} > 0$. In the second case, we need 
to have $z^*(\infty) = 0$ when $z \in M(M_{\phi}^*, F) \neq \{0\}$. Suppose, on the contrary, $z^*(\infty) = a > 0$. 
Since $z^* \in M(M_{\phi}^*, F)$ it follows that $\frac{a}{\phi} \leq \frac{z^*}{\phi} \in F$ and $1/\phi \in F$ gives, by maximality of the Marcinkiewicz space, that $1/\phi \in M_{\psi}^*$, where fundamental function of $F$ is $f_F(t) = \psi (t)$. It means that $\sup_{t > 0} \frac{\psi(t)}{\phi(t)} < \infty$. On the other hand, since $p_{\psi} \geq \alpha_F > q_{\phi}$ and $p_{\psi/\phi} \geq p_{\psi} - q_{\phi} > 0$ it follows that for $0 < \varepsilon < (p_{\psi} - q_{\phi} )/2$ and for large $t$ we obtain
\begin{equation*}
\frac{\psi(t)}{\phi(t)} \geq \frac{ \psi(1)}{\phi(1) \,m_{\phi}(t) m_{\psi}(1/t)} \geq \frac { \psi(1)}{\phi(1)} \, t^{p_{\psi}-\varepsilon - (q_{\phi} +\varepsilon)} = \frac{ \psi(1)}{\phi(1)} \, t^{p_{\psi} - q_{\phi} - 2 \varepsilon} \rightarrow \infty ~ {\rm as} ~ t \rightarrow \infty,
\end{equation*} 
a contradiction.
\end{proof}

\noindent
{\it Proof of Theorem 11.} (a) We have
\vspace{-2mm}
\begin{eqnarray*}
M_{\phi}^* \odot M(M_{\phi}^*, F) 
& \equiv&
L^{\infty}(\phi)^{(*)} \odot M[L^{\infty}(\phi)^{(*)}, F]  \hspace{11mm} ({\rm using ~ Lemma ~ 6}) \\
&=&
L^{\infty}(\phi)^{(*)} \odot F(1/\phi)^{(*)} \hspace{22mm} ({\rm by ~Theorem ~ 1(iv)}) \\
&\equiv&
\left\{ [L^{\infty}(\phi)^{(*)}]^{1/2} [F(1/\phi)^{(*)}]^{1/2}\right\}^{(1/2)} \hspace{2mm} ({\rm using ~ Lemma ~ 4}) \\
&=&
\left\{ [L^{\infty}(\phi)^{1/2} F(1/\phi)^{1/2}]^{(*)} \right\}^{(1/2)} \\
& &{\rm \hspace{1mm} (using ~the ~Krugljak-Maligranda ~result ~ [KM03], Thm ~ 2}) \\
&=&
\left\{ [(L^{\infty})^{1/2} F^{1/2}]^{(*)} \right\}^{(1/2)} \equiv \left\{ [F^{(2)}]^{(*)} \right\}^{(1/2)} \equiv F^{(*)} \equiv F.
\end{eqnarray*}

\vspace{-2mm}

\noindent
Note that Lemma 4 can be used in the above equality since  $0 < p_{\phi, I}  \leq q_{\phi, I} < 1$ implies 
that the operators $H, H^*: L^{\infty}(\phi) \rightarrow L^{\infty}(\phi)$ are bounded. Moreover, $\alpha_F > q_{\phi, I}$ gives that $H^*: F(1/\phi) \rightarrow F(1/\phi)$ is bounded and $\beta_ F < 1 +  p_{\phi, I}$, which is satisfied, that $H: F(1/\phi) \rightarrow F(1/\phi)$ is bounded (see [Ma80, Theorem 1] or [Ma83, Theorem 1).
\medskip

(b) By the duality results
\begin{eqnarray*}
E \odot M(E, \Lambda_{\phi}) 
&\equiv&
E \odot M((\Lambda_{\phi})^{\prime}, E^{\prime}) \equiv E \odot M(M_{t/\phi(t)}, E^{\prime}) \\
& \equiv&
E \odot M[L^{\infty}(\frac{t}{\phi(t)})^{(*)}, E^{\prime}]  \hspace{3mm} ({\rm using ~ Lemma ~ 6 ~and ~ symmetry ~ of} ~E) \\
&=&
E \odot E^{\prime}(\frac{\phi(t)}{t})^{(*)} \equiv E^{(*)} \odot E^{\prime}(\frac{\phi(t)}{t})^{(*)}  \hspace{6mm} ({\rm by ~Theorem ~ 1(iv)}) \\
&\equiv&
\left\{ [E^{(*)}]^{1/2} [E^{\prime}(\frac{\phi(t)}{t})^{(*)}]^{1/2}\right\}^{(1/2)} \hspace{2mm} ({\rm using ~ Lemma ~ 4}) \\
&=&
\left\{ [E^{1/2} E^{\prime}(\frac{\phi(t)}{t})^{1/2}]^{(*)} \right\}^{(1/2)} \\
&&{\rm \hspace{5mm} (using ~the ~Krugljak-Maligranda ~result ~ [KM03], Thm ~ 2}) \\
&=&
\left\{ [E^{1/2} (E^{\prime})^{1/2} (t^{-1/2} \phi(t)^{1/2})]^{(*)} \right\}^{(1/2)} \hspace{6mm} ({\rm by ~Theorem ~ 1(iv)}) \\
&\equiv&
\left[ (E \odot E^{\prime})^{(2)} (t^{-1/2} \phi(t)^{1/2})^{(1/2)} \right]^{(*)} \\
& & \hspace{5mm} ({\rm by ~the ~Lozanovskii ~ factorization ~theorem)}\\
&\equiv&
\left[ L^2(t^{-1/2} \phi(t)^{1/2})^{(1/2)} \right]^{(*)} \equiv L^1(\frac{\phi(t)}{t})^{(*)} \equiv \Lambda_{\phi, 1}.
\end{eqnarray*}
We must only control if the assumptions from Lemma 4 are satisfied in our case, that is, if operators $H, H^*$ are bounded in $E$ and in $E^{\prime}(\frac{\phi(t)}{t})$. Since $0 < \alpha_E \leq \beta_E < p_{\phi, I} < 1$ it follows from Boyd's result that $H$ and $H^*$ are bounded in $E$ (cf. [KPS82], pp. 138-139, [BS88], Theorem 5.15 and [KMP07], pp. 126-129). The boundedness of $H$ in $E^{\prime}(\frac{\phi(t)}{t})$ is equivalent to the estimate 
$$
\| \frac{\phi(t)}{t^2} \int_0^t x(s) \frac{s}{\phi(s)}\, ds \|_{E^{\prime}} \leq C_1 \, \| x \|_{E^{\prime}} ~ {\rm for ~all}  ~ x \in E^{\prime},
$$
which is true if $\beta_{E^{\prime}} < p_{t^2/\phi(t), I}$ (cf. [Ma80, Theorem 1] or [Ma83, Theorem 1). The last strict inequality means that $\beta_{E^{\prime}} = 1 - \alpha_E < p_{t^2/\phi(t), I} = 2 - q_{\phi, I}$ or $\alpha_E > q_{\phi, I} - 1$ which is true because $\alpha_E > 0$ and $q_{\phi, I} < 1$. The boundedness of $H^*$ in $E^{\prime}(\frac{\phi(t)}{t})$ is equivalent to the estimate $$
\| \frac{\phi(t)}{t} \int_t^l x(s) \frac{1}{\phi(s)}\, ds \|_{E^{\prime}} \leq C_2 \, \| x \|_{E^{\prime}} ~ {\rm for ~all}  ~ x \in E^{\prime},
$$
which is true if $\alpha_{E^{\prime}} > 1- p_{\phi, I}$ (cf. [Ma80, Theorem 1] or [Ma83, Theorem 1]). The last strict inequality means that $\alpha_{E^{\prime}} = 1 - \beta_E > 1- p_{\phi, I}$ or $\beta_E < p_{\phi, I}$, but this is true by the assumption.
\qed

\vspace{3mm}
{\bf Examples 5.} (a) If $E = L^q, F = L^p$ and $\phi(t) = t^{1/r}$, where $1 \leq p < r < q < \infty$, then from Theorem 11 we obtain 
\begin{equation} \label{equation38}
L^{r, \infty} \odot M(L^{r, \infty}, L^p) = L^p ~ {\rm and } ~~ L^q \odot M(L^q, L^{r, 1}) = L^{r, 1}.
\end{equation}

\vspace{-2mm}

\noindent
Equalities (\ref{equation38}) we can also get from  (\ref{equation31}) and (\ref{equation30}). In fact, space $L^{r, \infty}$ 
satisfies upper $r$-estimate (cf. [Ma04], Theorem 5.4(a) and [KK05], Theorem 3.1 and Corollary 3.9), therefore  for 
$p < r$ is $p$-convex with some constant $K \geq 1$ (cf. [LT79], Theorem 1.f.7). After renorming it is $p$-convex with 
constant $1$ and we are getting from (\ref{equation31}) the first equality in (\ref{equation38}) with equivalent norms. 
On the other hand, $L^{r,1}$ satisfies lower $r$-estimate (cf. [Ma04], Theorem 5.1(a)), therefore for $q > r$ it is $q$-concave 
with some constant $K \geq 1$ (cf. [LT79], Theorem 1.f.7). After renorming is $q$-concave with constant $1$ and we are 
getting from (\ref{equation30}) the second equality in (\ref{equation38}) with equivalent norms. 

(b)  If $E = L^{q, r}, F = L^{p, r}$ and $\phi(t) = t^{1/s}$, where $1 \leq p < s < q < \infty$ and $1 \leq r \leq \infty$, then from Theorem 11 we obtain 
\vspace{-2mm}
\begin{equation} \label{equation39}
L^{s, \infty} \odot M(L^{s, \infty}, L^{p, r}) = L^{p, r} ~ {\rm and } ~~ L^{q, r} \odot M(L^{q, r}, L^{s, 1}) = L^{s, 1}.
\end{equation}

(c) If $F = \Lambda_{\psi}$ and $\alpha_F = p_{\psi, I} > q_{\phi, I}$, then from Theorem 11(a) we also obtain factorization 
(\ref{equation34}) since $p_{\psi/\phi, I} \geq p_{\psi, I} - q_{\phi, I} > 0$.

\vspace{2mm}
{\bf Remark 8.} In the case $I = (0, 1)$ the assumption $\alpha_F > q_{\phi, I}$ implies the imbedding $M_{\phi}^* \hookrightarrow F$, even the imbedding $M_{\phi}^* \hookrightarrow \Lambda_{\psi}$, where $\psi$ is a fundamental function of $F$ because $p_{\psi, 1} \geq \alpha_F > q_{\phi, 1}$ and $p_{\psi/\phi, 1} \geq p_{\psi, 1} - q_{\phi, 1} > 0$ gives
\vspace{-2mm}
$$
\int_0^1 \frac{1}{\phi(t)} \psi^{\prime}(t)\, dt \leq \int_0^1 \frac{\psi(t)}{\phi(t)} \, \frac{dt}{t} < \infty.
$$
Consequently, $M(M_{\phi}^*, F) \neq \{0\}$.


\vspace{3mm}

\noindent {\footnotesize Pawe\l\ Kolwicz and Karol Le\'{s}nik, Institute of
Mathematics of Electric Faculty\newline
Pozna\'n University of Technology, ul. Piotrowo 3a, 60-965 Pozna\'{n}, Poland
}\newline
\textit{E-mails:} ~\texttt{pawel.kolwicz@put.poznan.pl,
klesnik@vp.pl}\newline

\vspace{-1mm}

\noindent {\footnotesize Lech Maligranda, Department of Engineering Sciences and Mathematics\\
Lule\aa\ University of Technology, SE-971 87 Lule\aa , Sweden}\newline ~\textit{E-mail:} \texttt{lech.maligranda@ltu.se
}

\end{document}